\pdfoutput=1
\RequirePackage{ifpdf}
\ifpdf % We~are running pdfTeX in pdf mode
\documentclass[pdftex]{sigma}
\else
\documentclass{sigma}
\fi

\usepackage{enumerate}

\numberwithin{equation}{section}

\newtheorem{Theorem}{Theorem}[section]
\newtheorem{Lemma}[Theorem]{Lemma}
\newtheorem{Fact}[Theorem]{Fact}
\newtheorem{Cor}[Theorem]{Corollary}
\newtheorem{Prop}[Theorem]{Proposition}
\newtheorem{State}[Theorem]{Statement}

 { \theoremstyle{definition}
\newtheorem{Definition}[Theorem]{Definition}
\newtheorem{Remark}[Theorem]{Remark} }

\def\II{{\rm II}}
\def\tr{{\rm tr}}

\def\R{\mathbb{R}}
\def\B{\mathbb{B}}
\def\S{\mathbb{S}}
\def\H{\mathbb{H}}
\def\T{\mathbb{T}}
\def\Z{\mathbb{Z}}

\def\D{\mathbb{D}}

\def\Bb{\mathbf{B}}

\def\bs{\boldsymbol}

\def\refOmega{\Omega_0}

\def\Ucal{\mathcal{U}}

\def\Hcal{\mathcal{H}}
\def\Acal{\mathcal{A}}
\def\Ccal{\mathcal{C}}
\def\Dcal{\mathcal{D}}

\def\Bcal{\mathcal{B}}
\def\Scal{\mathcal{S}}
\def\Ecal{\mathcal{E}}

\def\Ncal{\mathcal{N}}

\def\Ocal{\mathcal{O}}

\def\Kcal{\mathcal{K}}

\def\supp{\operatorname{supp}}

\def\ed{{\rm d}}

\def\Ric{{\rm Ric}}
\def\area{{\rm area}}

\def\id{{\rm id}}
\def\Rm{{\rm Rm}}
\def\tr{{\rm tr}}

\def\dist{{\rm dist}}

\def\Rm{{\rm Rm}}
\def\tilRm{\widetilde\Rm}

\def\<{\langle}
\def\>{\rangle}

\def\setdiff{\backslash}

\def\Nsep{\textbf{\textup{NSep}$^+$}}

\begin{document}
\allowdisplaybreaks

\renewcommand{\thefootnote}{}

\newcommand{\arXivNumber}{2303.15752}

\renewcommand{\PaperNumber}{083}

\FirstPageHeading

\ShortArticleName{Rigidity and Non-Rigidity of $\mathbb{H}^n/\mathbb{Z}^{n-2}$ with Scalar Curvature Bounded from Below}

\ArticleName{Rigidity and Non-Rigidity of $\boldsymbol{\mathbb{H}^n/\mathbb{Z}^{n-2}}$\\ with Scalar Curvature Bounded from Below\footnote{This paper is a~contribution to the Special Issue on Differential Geometry Inspired by Mathematical Physics in honor of Jean-Pierre Bourguignon for his 75th birthday. The~full collection is available at \href{https://www.emis.de/journals/SIGMA/Bourguignon.html}{https://www.emis.de/journals/SIGMA/Bourguignon.html}}}

\Author{Tianze HAO~$^{\rm a}$, Yuhao HU~$^{\rm ab}$, Peng LIU~$^{\rm a}$ and Yuguang SHI~$^{\rm a}$}

\AuthorNameForHeading{T.~Hao, Y.~Hu, P.~Liu and Y.~Shi}

\Address{$^{\rm a)}$~Key Laboratory of Pure and Applied Mathematics,
School of Mathematical Sciences,\\
\hphantom{$^{\rm a)}$}~Peking University, Beijing, 100871, P.R.~China}
\EmailD{\href{mailto:haotz@pku.edu.cn}{haotz@pku.edu.cn},
		\href{mailto:1801110011@pku.edu.cn}{1801110011@pku.edu.cn},
		\href{mailto:ygshi@math.pku.edu.cn}{ygshi@math.pku.edu.cn}}

\Address{$^{\rm b)}$~School of Mathematical Sciences, Shanghai Jiao Tong University,\\
\hphantom{$^{\rm b)}$}~Shanghai, 200240, P.R.~China}
\EmailD{\href{mailto:yuhao.hu@sjtu.edu.cn}{yuhao.hu@sjtu.edu.cn}}

\ArticleDates{Received April 08, 2023, in final form October 20, 2023; Published online November 01, 2023}

\Abstract{We show that the hyperbolic manifold $\mathbb{H}^n/\mathbb{Z}^{n-2}$ is not rigid under all compactly supported deformations that preserve the scalar curvature lower bound $-n(n-1)$, and that it is rigid under deformations that are further constrained by certain topological conditions. In addition, we prove two related splitting results.}


\Classification{53C21; 53C24}

\begin{flushright}
\begin{minipage}{65mm}
\it Dedicated to Jean-Pierre Bourguignon\\ on the occasion of his 75th birthday
\end{minipage}
\end{flushright}

\renewcommand{\thefootnote}{\arabic{footnote}}
\setcounter{footnote}{0}

\section{Introduction}

In \cite[Section 3]{Gr2019} and \cite[p.~240]{Gr2021}, Gromov stated the following generalization
of Min-Oo's hyperbolic rigidity theorem \cite{MinOo89}.

\begin{State}[``generalised Min-Oo rigidity theorem'']\label{generminoo}
Parabolic quotients $Z=\mathbb{H}^n/\Gamma$ of the hyperbolic $n$-space admit no non-trivial, compactly supported `deformation' with scalar curvature $R\geq -n(n-1)$.
\end{State}
According to \cite{Gr2019}, a \emph{deformation} can change not only the metric, but also the topology of a compact region in $Z$.
If the deformation is topologically a connected sum with
a closed $n$-manifold,
Statement~\ref{generminoo} is known to be true for (at least) $Z = \H^n/\Z^{n-1}$,
with idea of proof already outlined by \cite[Section $5\frac{5}{6}$]{Gr93} (for a detailed treatment, see
also  \cite[Theorem 1.1]{ACG08}).
The situation
turns out to be more subtle
if broader types of deformations are considered, allowing, for example, surgeries along
an embedded, non-contractible loop.
In this latter case we construct a counterexample to Statement~\ref{generminoo}, which, more precisely,
demonstrates the following.

\begin{Theorem}\label{counterexample}
For $n\ge 3$, let $\mathbb{H}^n/\mathbb{Z}^{n-2}$ be equipped with the standard hyperbolic metric.
There exists a complete Riemannian manifold $(M^n,g)$, not $($globally$)$ hyperbolic, and
compact subsets $K\subset M$ and $K'\subset\H^n/\Z^{n-2}$, such that $(1)$~$R_g\ge -n(n-1)$ and $(2)$~$M\setminus K$ is isometric to~$\big(\mathbb{H}^n/\mathbb{Z}^{n-2}\big)\setminus K'$.	
\end{Theorem}
\begin{Remark}\quad
\begin{enumerate}\itemsep=0pt
	\item[$(1)$] While the theorem above concerns non-rigidity of $\H^n/\Z^{n-2}$, it is also interesting to
    ask whether its statement still holds if $\H^n/\Z^{n-2}$ is replaced by $\H^n/\Z^{n-1}$; this will be answered in the affirmative in Section~\ref{CErmksurg}.
    Thus, we obtain counterexamples to the ``weak rigidity of $\mathbb{H}^n/\mathbb{Z}^{n-1}$'' mentioned in \cite[p.~678]{Gromov18GAFA}.

   	 \item[$(2)$] En route to proving Theorem \ref{counterexample}, we obtain counterexamples
    (see Proposition \ref{XrbarInfo}) to the following statement in \cite[p.~12]{Gr2019}: \emph{Represent $\H^n/\Z^{n-2}$ as a warped product $\bigl(\H^2\times \T^{n-2}, g_H\bigr)$ $($see formula~\eqref{gHdef}$)$, and, for a geodesic $2$-disk $\D^2\subset\H^2$, let $X = \D^2 \times \T^{n-2}\subset \H^2\times \T^{n-2}$ with the restricted metric $g_H|_X$; then no Riemannian manifold $\big(M^n,g\big)$ with boundary isometric to $\partial X$ can have scalar curvature $R_g\ge -n(n-1)$ and mean curvature\footnote{Unless specified otherwise, in this article the mean curvature along a boundary will always be computed with respect to the \emph{outward} unit normal.} of $\partial M$ greater than that of $\partial X$.}

    \item[$(3)$] Our proof of Theorem~\ref{counterexample} is constructive, which, a little to our surprise, shows that $M$ can be chosen to be homeomorphic to $\H^n/\Z^{n-2}$ (see Section~\ref{CErmktop}); moreover, $R_g> -n(n-1)$ for some
    points in $K$.
\end{enumerate}
\end{Remark}

From the perspective of our construction, the non-rigidity of $\H^n/\Z^{n-2}$ seems closely related to the fact:
\emph{A deformation supported in a compact subset $K$ can
`break' the incompressibility\,\footnote{A continuous map $f\colon X\rightarrow Y$ between topological spaces
is said to be \emph{incompressible} if the induced map $f_*\colon \pi_1(X)\rightarrow \pi_1(Y)$ is injective;
when $f$ is an inclusion, we say \emph{`$X$ is incompressible in $Y$'}.} of some submanifold that is disjoint from $K$.}
On the other hand, rigidity does hold if one only considers deformations
that preserve such incompressibility, as the next theorem shows (cf. \cite[Theorem 1.8]{CLSZ2021}).

\begin{Theorem}\label{deform1}
For $3\leq n\leq 7$, let $(M^n,g)$ be a complete Riemannian manifold\,\footnote{In this article,
all manifolds are assumed to be orientable, and all hypersurfaces $2$-sided.}
with  scalar curvature $R_g\geq-n(n-1)$. Suppose that there exist compact subsets $K\subset M$, $K'\subset \H^n/\Z^{n-2}$, and an isometry $f\colon M\setminus K\rightarrow \bigl(\mathbb{H}^n/\mathbb{Z}^{n-2}\bigr)\setminus K'$. Representing $\H^{n}/\Z^{n-2}$ topologically as $\R^2_+\times \T^{n-2}$, let $p\in \R^2_+$ be such that
$T = \{p\}\times \T^{n-2}$ is disjoint from $K'$, and suppose that the map $f^{-1}|_T\colon T\rightarrow M$ is incompressible.
Then $(M,g)$ is isometric to $\mathbb{H}^n/\mathbb{Z}^{n-2}$.
\end{Theorem}

Technically, we will derive Theorem~\ref{deform1}
as a consequence of Theorem~\ref{alhpmt4} below. The latter can be regarded as a kind of positive mass type theorem for manifolds with an ALH end; its
 statement relies on a gluing construction, which we now describe.

\emph{Gluing construction:} Let $N^n$ be a smooth manifold, and
suppose that $\phi\colon \T^k\rightarrow N$ $(1\le k\le n-2)$ is an embedding with trivial normal bundle.
Moreover, write
$\H^n/\Z^{n-1}$ (topologically) as the product $\R\times \T^{n-k-1}\times \T^k$,
and define
\[
	\psi\colon \ \T^k\rightarrow \R\times \T^{n-k-1}\times \T^k  \cong \H^n/\Z^{n-1}
	\qquad
	\mbox{by}
	\quad \psi(p) = (t,q,p)
\]
for some fixed $t\in \R$ and $q\in \T^{n-k-1}$.
By removing tubular neighborhoods of $\phi\bigl(\T^k\bigr)\subset N$ and $\psi\bigl(\T^k\bigr)\subset \H^{n}/\Z^{n-1}$
and then identifying the respective boundaries in the obvious way, we obtain
a manifold $M$. For brevity, $M$ will be referred to as obtained by \emph{gluing
$N$ and $\H^n/\Z^{n-1}$ along~$\T^k$ via~$(\phi,\psi)$.} In particular, for $c$ sufficiently large,
$(c,\infty)\times \T^{n-1}\subset \H^n/\Z^{n-1}$ remains an `end' of~$M$,
and this end is denoted by $\Ecal$.

\begin{Theorem}\label{alhpmt4}
For $3\leq n\leq 7$, let $N^n$ be a smooth manifold that is either closed or non-compact without boundary, and let
$M^n$ be obtained by gluing $N$ with $\H^n/\Z^{n-1}$ along $\T^{k}$ via $(\phi,\psi)$  $($see description above$)$.
Suppose that
\begin{enumerate}[$(a)$]\itemsep=0pt
	\item the map $\phi\colon \T^k\rightarrow N$ is incompressible; \label{incprsAssu}
\item $g$ is a complete Riemannian metric on $M$ with $R_g\geq -n(n-1)$;%\label{scalarbdAssu}
\item  $(\Ecal, g)$ is asymptotically locally hyperbolic $($ALH\,$)$ $($see Definition~{\rm \ref{alh})}.	
	%\label{ALHAssu}
\end{enumerate}
Then $\bar m_{\Ecal,g}\ge 0$ $($see Definition~{\rm \ref{massDef})}. In addition, suppose that
\begin{enumerate}[$(a)$]\itemsep=0pt\setcounter{enumi}{3}
\item the curvature tensor of $(M,g)$ and its first covariant derivatives are bounded;\label{curvbdAssu}
\item there exists some $\alpha>0$ such that $R_g \leq -\alpha$ outside a compact set.
	\label{RbdAssu}
\end{enumerate}
Then $\kappa = 0$ $($see \eqref{alhasymp} for the definition of $\kappa$$)$ only if $(M,g)$ is Einstein.
\end{Theorem}

Readers familiar with positive mass theorems may
have noticed that the second half of Theorem~\ref{alhpmt4}
is not in an ideal form; in other words, one wants to know whether
the vanishing of~$\bar m_{\Ecal,g}$, and not just $\kappa$, implies that $(M,g)$ is isometric to $\H^n/\Z^{n-1}$, even without the assumptions~$(d)$ and~$(e)$. In our proof, these assumptions play a role in making sure that the
normalized Ricci flow (NRF) starting at $g$ has desired properties (see Lemma~\ref{nrfdeform});
on the other hand,
it seems subtle to prove
hyperbolicity from $(M,g)$ being Einstein and the assumed ALH decay rate. Thus we decide to~leave the stronger statement for future investigation.

Theorem~\ref{alhpmt4} has the following corollary.

\begin{Cor}\label{hpbInfCor}
For $3\leq n\leq 7$, let $N^n$ be a closed manifold, and suppose that
$M^n$ is obtained by gluing $N$ with $\H^n/\Z^{n-1}$ along $\T^{k}$ via $(\phi,\psi)$.
Suppose that $g$ is a complete metric on $M$ such that $(M,g)$ is isometric to the hyperbolic manifold $\H^{n}/\Z^{n-1}$ outside a
 compact set,\footnote{That is, there exists an isometry $f\colon M\setminus K\rightarrow(\H^n/\Z^{n-1})\setminus K'$ for some compact
 sets $K\subset M$ and $K'\subset \H^n/\Z^{n-1}$.} and suppose that
 \begin{enumerate}\itemsep=0pt
\item[$(a)$]  the map $\phi\colon \T^k\rightarrow N$ is incompressible;
\item[$(b)$]  $R_g\geq -n(n-1)$.
\end{enumerate}
Then $(M,g)$ is isometric to $\H^n/\Z^{n-1}$.
\end{Cor}

In fact, Corollary~\ref{hpbInfCor} remains true if $N$ is allowed to be non-compact,
which can be deduced as a corollary of Theorem~\ref{alhpmt1} below (see Remark~\ref{hpbInfCor_NonCptRmk}).

Besides rigidity problems modeled on complete manifolds, it is often natural to consider similar problems for manifolds with boundary and scalar/mean curvature bounds.
In this regard, we present a splitting result of `cuspidal-boundary' type
(see \cite[Section 4, last paragraph]{Gr2019}). Our proof relies on an approximation scheme developed in \cite{Zhu2020} involving
$\mu$-bubbles.

\begin{Theorem}\label{alhpmt2}
Let $\bigl(M^4, g\bigr)$ be a complete, non-compact Riemannian $4$-manifold with
compact, connected boundary $\partial M$.
Suppose that $\pi_2(M) = \pi_3(M) = 0$ and that the scalar curvature $R_g\geq -12$. Then
\[
\inf_{\partial M}H \leq 3,
\]
where $H$ is the mean curvature of $\partial M$.
Moreover, if
\[
\inf_{\partial M}H=3,
\]
then $(M, g)$ is isometric to $\big((-\infty, 0]\times \Sigma, \ed t^2+{\rm e}^{2t}g_0\big)$, where $t$ is the coordinate on $(-\infty, 0]$ and $(\Sigma, g_0)$ is a closed
$3$-manifold with a flat metric.
\end{Theorem}

\begin{Remark}
	Theorem~\ref{alhpmt2} would fail if one allows $M$ to be compact. Indeed,
	take
\[
M=\S^1\times\B^3, \qquad
g=\cosh^2\rho  \ed\theta^2 +\ed\rho^2 +\sinh^2\rho  g_{\S^2},\qquad \rho\leq \rho_0,
\]
where $\theta\in \S^1$, $\rho$ is the radial coordinate on $\B^3$, and $g_{\S^2}$ is the standard round metric on $\S^{2}$. In~this example, $M$ has contractible universal cover, so both
its $\pi_2$ and $\pi_3$ vanish. Moreover, since~$g$ is hyperbolic, $R_g = -12$,
but the mean curvature $H_{\partial M} = 2\coth\rho_0+\tanh\rho_0>3$.

Counterexamples also exist if one drops the assumption on $\pi_2(M)$ and $\pi_3(M)$. In fact, let us take the manifold
$(\bs{M}', \bs{g}')$ in Section~\ref{CErmksurg} and then, for sufficiently small $z_0>0$, remove the subset $\{0<z<z_0\}$ from $\bs{M}'$; the
result is a manifold $\bs{M}''$ with \[ \pi_2(\bs{M}'')\ne 0, \qquad H_{\partial\bs{M}''} = 3\qquad \mbox{and} \qquad R \ge -12.\]
Clearly, $\bs{M}''\not\cong [c,\infty)\times \partial{\bs{M}''} \cong [c,\infty)\times \T^3$.
\end{Remark}

Finally, we present an analogue of Theorem \ref{alhpmt2} in more general dimensions.
\begin{Definition}[{cf.~\cite{CLSZ2021}}]\label{cdegdef}
We say that a closed, connected manifold $\Sigma$ \emph{belongs to the class $\mathcal{C}_{\deg}$},~if
\begin{itemize}\itemsep=0pt
	\item $\Sigma$ is aspherical,\footnote{A closed, connected manifold is said to be \emph{aspherical} if
it has contractible universal cover.} and
	\item any compact manifold $\Sigma'$ that admits a map to $\Sigma$ of nonzero degree
			cannot be endowed with a PSC metric (i.e., metric with positive scalar curvature).
\end{itemize}
\end{Definition}
It is well known that $\T^n\in \Ccal_{\deg}$ for $n\le 7$; also note that the second item
in Definition~\ref{cdegdef} is redundant when $\dim \Sigma \le 5$,
according to \cite{CLL2021}.

\begin{Theorem}\label{alhpmt1}
For $3\leq n\leq 7$, let $(M^n, g)$ be a complete and non-compact Riemannian manifold with compact, connected boundary $\partial M$. Suppose that
\begin{enumerate}\itemsep=0pt
\item[$(a)$] $\partial M$ is incompressible in M,
\item[$(b)$] $\partial M\in\Ccal_{\deg}$,
\item[$(c)$]  $R_g\geq -n(n-1)$,
\end{enumerate} then
\begin{equation*}%\label{meancurvineq1}
\inf_{\partial M}H \leq n-1,
\end{equation*}
where $H$ is the the mean curvature of $\partial M$.

Moreover, if
\[
\inf_{\partial M}H = n-1,
\]
then $(M, g)$ is isometric to $((-\infty, 0]\times \Sigma, \ed t^2+{\rm e}^{2t}g_0)$ where $t$ is the coordinate on $(-\infty, 0]$ and $(\Sigma, g_0)$ is a closed
$(n-1)$-manifold with a flat metric.
\end{Theorem}

\noindent\textbf{Additional notes on the literature.}
$(a)$ All our main theorems are fundamentally related to Gromov's fill-in problems
	(e.g., \cite[Problems A and B]{Gr2019}; \cite[p.~234, Question~(c)]{Gr2021}).
$(b)$~Theorem~\ref{alhpmt1} can be viewed as a generalization of \cite[Theorem~3.2]{Yau2001}.
$(c)$ It is a classical theme to relate incompressibility conditions with scalar curvature (see \cite[Section~11]{GL1983}).
$(d)$~To adapt to modern language, our Theorem~\ref{alhpmt4} considers manifolds with a prescribed end
and some `arbitrary ends'; the study of positive-mass type theorems on such manifolds has generated
considerable interest recently (see, for example, \cite{CZ21pmt,CL20,LUY21,Zhu22PMT}).
$(e)$~While in this paper we focus on rigidity results for complete, non-compact
manifolds with boundary and scalar curvature lower bounds, similar results in the
compact case (with boundary) are obtained by Gromov in \cite[Section~4]{Gr2019}.
In both cases, the proofs rely on the $\mu$-bubble technique.
$(f)$ It would be interesting to compare Theorem~\ref{alhpmt4}
with some recent progress in proving positive mass and rigidity results for
ALH manifolds (see \cite{AHK2022,CG2021,CGNP2018,HJ2022});
in this latter development, manifolds are often assumed to have nonempty inner boundary with the
mean curvature bound $H\le n-1$ (now~$H$ is computed with respect to the \emph{inner} unit normal);
 such mean curvature bounds serve as barrier conditions in the method of `marginally outer
trapped surfaces' (MOTS), which can be viewed as a generalization of the $\mu$-bubble technique.

\textbf{Organization of this article.}
The proof of Theorem~\ref{counterexample} is technically independent from
the~rest of the work and is included in Section~\ref{Non-rigidity}. Section~\ref{ALH_manifolds} serves as a preliminary to proving
Theorem~\ref{alhpmt4}, presenting results concerning NRF and conformal deformations.
In Section~\ref{Two_rigibility}, we prove Theorem~\ref{alhpmt4}, followed by
proofs of Corollary~\ref{hpbInfCor} and Theorem~\ref{deform1}.
In Section~\ref{Two_splitting}, we prove Theorem~\ref{alhpmt2} and Theorem~\ref{alhpmt1}.
Several of the proofs rely on the so-called `$\mu$-bubble' technique,
a brief discussion of which is included in Appendix~\ref{mubbSec}. Appendix~\ref{TopologicalLemmas} includes
two topological lemmas.

\section[Non-rigidity of H\^{}n/Z\^{}(n-2)]{Non-rigidity of $\boldsymbol{\H^n/\Z^{n-2}}$}\label{Non-rigidity}

Let the hyperbolic $n$-space $\H^n$ be represented by the upper half-space model
$\R^n_+ = \{(x,\boldsymbol{y},z)$: $x\in \R, \boldsymbol{y}\in \R^{n-2}, z>0\}$,
and let $\Z^{n-2}$ act by translating along the orthogonal lattice $2\pi \Z^{n-2}\subset \R^{n-2}$
while keeping the $x,z$-coordinates fixed. The quotient space is denoted by
$\H^{n}/\Z^{n-2}$ and has the hyperbolic metric
\begin{equation}\label{gHdef}
	g_{H} = z^{-2}\bigl(\ed z^2 + \ed x^2\bigr) + z^{-2} g_{\T^{n-2}},
\end{equation}
where the subscript `$H$' stands for `hyperbolic', and $g_{\T^{n-2}}$ is the associated flat metric on~$\T^{n-2}$.
Henceforth, we will regard $(x,z)$ as coordinates on the hyperbolic plane $\H^2$;
manifestly that $\big(\H^{n}/\Z^{n-2},g_H\big)$ is a warped product of $\H^2$ and $\big(\T^{n-2}, g_{\T^{n-2}}\big)$.

The following lemma is easily verified by standard computation, so we omit its proof.
\begin{Lemma}\label{H2DerHess}
Let $\nabla$, $\nabla^2$ denote the gradient and Hessian
with respect to $g_H$ (same below). We~have
\begin{enumerate}\itemsep=0pt
    \item[$(a)$] $\nabla z=z^2{\partial}/{\partial z}$,
	\item[$(b)$] $\nabla^2z({\partial}/{\partial x}, {\partial}/{\partial x}) = -\nabla^2z({\partial}/{\partial z}, {\partial}/{\partial z})=-1/z$,
    \item[$(c)$] $\nabla^2z({\partial}/{\partial z}, {\partial}/{\partial x})=0$.
\end{enumerate}	
\end{Lemma}

Next, we proceed to prove
Theorem~\ref{counterexample} by constructing an example
that satisfies all its conditions. The idea is to remove a suitable compact
subset, $X_{p,r}$, from $\H^n/\Z^{n-2}$ and then `glue' the result with a
compact manifold, $\bar X_r$, along their boundaries; $X_{p,r}$ and
$\bar X_r$ will be defined in Sections~\ref{1stprelimSec} and \ref{2ndprelimSec}
respectively, and then we handle the gluing step
in Section~\ref{gluingSec}.

\subsection{First preliminary construction}\label{1stprelimSec}
Let $p\in \H^2$, and define
\begin{equation}\label{Xrp_def}
	X_{p,r}:= \D_r(p)\times \T^{n-2}\subset \H^{n}/\Z^{n-2} \qquad
	\mbox{and} \qquad  Y_{p,r} :=\partial X_{p,r},
\end{equation}
where $\D_r(p)\subset \H^2$ is the geodesic disc, centered at $p$, of radius $r>0$;
the inclusion in \eqref{Xrp_def} makes sense since $\H^n/\Z^{n-2}$ is a warped product
of $\H^2$ and $\T^{n-2}$, as we already noted.

Now we have two sets of coordinates for $\H^2$: $(x,z)$ and the polar coordinates
$(\varrho,\theta)$ centered at $p$. In terms of the polar coordinates, the metric on $\H^2$ reads
\begin{equation*}%\label{H2metricpolar}
	g_{\H^2} = \ed\varrho^2 + \sinh^2\varrho\,\ed\theta^2.
\end{equation*}

\begin{Lemma}\label{YprMean}
The boundary $Y_{p,r}\subset (X_{p,r}, g_H)$ has the mean curvature
\begin{equation}\label{Hpr}
	H_{p,r}=\coth r-(n-2)z^{-1}\frac{\partial z}{\partial \varrho}.
\end{equation}
Moreover,
\begin{enumerate}\itemsep=0pt
    \item[$(a)$] $|H_{p,r} - \coth r| \leq n-2$;
	\item[$(b)$] There exists a constant $r_0>0$ such that $H_{p,r}>0$ for all $r\leq r_0$.
\end{enumerate}

\end{Lemma}
\begin{proof}
The formula \eqref{Hpr} is straightforward to check by using the representation
\[
	g_H = \ed\varrho^2 + \sinh^2\varrho\,\ed\theta^2 + z^{-2}g_{\T^{n-2}}.
\]
Moreover, since both $z^{-1}\nabla z$ and $\nabla \varrho$ have unit norm with respect to $g_H$,
\begin{equation}\label{dzdrhoBd}
\left|\frac{\partial z}{\partial \varrho}\right|
= \left|\bigl\langle\nabla z, \nabla \varrho \bigr\rangle\right| =\left|z\bigl\langle z^{-1}\nabla z, \nabla \varrho \bigr\rangle\right|\le z.
\end{equation}
This implies $(a)$, and $(b)$ follows since $\coth r\rightarrow\infty$ as
$r\rightarrow 0$.
\end{proof}

\begin{Lemma}\label{dthetazEst}
	There exists a constant $C_r>0$, depending only on $r$, such that on $\partial\D_r(p)$ we have
	\begin{equation*}
		\left|\partial_\theta z (r,\theta)\right| \le z \sinh r \qquad
		\mbox{and}\qquad
		\left|\partial_\theta^2 z(r,\theta)\right| \le C_r z.
	\end{equation*}
\end{Lemma}
\begin{proof}
	Since both $z^{-1}\nabla z$ and $(\sinh r)^{-1}(\partial/\partial\theta)$
	have unit norm with respect to $g_{H}$, we have
	\[
		|\partial_\theta z(r,\theta)| = \left|\bigl\langle\nabla z, \partial/\partial\theta\bigr\rangle\right|
		\le z\sinh r.
	\]
	Moreover, a calculation shows that
	\begin{equation}\label{hessztt}
		\nabla^2 z (\partial/\partial\theta,\partial/\partial\theta)
		= \partial_\theta^2 z  +(\partial_\varrho z) \sinh \varrho \cosh \varrho.
	\end{equation}
	By Lemma~\ref{H2DerHess}\,$(b)$,\,$(c)$, the left-hand side of \eqref{hessztt} has its magnitude bounded
	by $(\sinh^2\rho)z$; thus, using~\eqref{dzdrhoBd} and evaluating
	\eqref{hessztt} at $\varrho = r$, we get
	\[
		\big|\partial_\theta^2 z(r,\theta)\big| \le \sinh r (\sinh r + \cosh r) z.
	\]
	Taking $C_r  = \sinh r(\sinh r + \cosh r)$ finishes the proof.
\end{proof}

\subsection{Second preliminary construction}\label{2ndprelimSec}
Let $\Dcal$ be a $2$-disc with polar coordinates $\bigl(\bar\varrho,\bar\theta\bigr)$, where
\[
	0\le \bar\varrho\le \pi/3 \qquad \mbox{and}\qquad 0\le \bar \theta<2\pi.
\]
Equip $\Dcal$ with the metric
\[
	g_\Dcal = \ed\bar\varrho^2 + 4\sin^2(\bar\varrho/2)\ed\bar\theta^2.
\]
Thus, $(\Dcal, g_\Dcal)$ is isometric to a `cap' in the round sphere of radius $2$.

Now let $r >0$ and $z(\varrho,\theta)$ be as in Section~\ref{1stprelimSec} above. Consider
\begin{equation*}
	\bar X_r := \S^1 \times \Dcal \times \T^{n-3}
\end{equation*}
equipped with the metric
\begin{equation}\label{barg}
\bar g =\sinh^ 2r \,\ed\theta^2 +\bigl(z(r,\theta)\bigr)^{-2} g_\Dcal+\bigl(z(r,\theta)\bigr)^{-2}g_{\T^{n-3}},
\end{equation}
and let $\bar Y_r := \partial \bar X_r$.
By construction, the boundaries $(Y_{p,r}, g_H|_{Y_{p,r}})$ and $\bigl(\bar Y_r, \bar g|_{\bar Y_r}\bigr)$ are isometric under the obvious identification.

\begin{Lemma}%\label{YrMean}	
The boundary $\bar Y_r\subset \bigl(\bar X_r, \bar g\bigr)$ has the mean curvature
\begin{equation}\label{barHr}
\bar H_r=\frac{\sqrt{3}}{2}z(r,\theta).
\end{equation}
\end{Lemma}
\begin{proof}
Standard computation by using \eqref{barg}.
\end{proof}

Regarding the scalar curvature of a warped-product metric,
the following is well-known.

\begin{Lemma}[{cf.~\cite[Proposition 7.33]{GL1983}}]\label{scalarwarp}
Let $\bigl(N^{n-1}, h\bigr)$ be an $(n-1)$-dimensional Riemannian manifold with scalar curvature $R_h$. Given any smooth function $\phi(\theta)$ defined on an
interval $I$ and a constant $a>0$, the warped product metric
$
	g = a^2\ed\theta^2+ \phi(\theta)^2 h
$
defined on $I\times N$ has the scalar curvature
\begin{equation}\label{scalarWarpEq}
R_g=\frac{n-1}{a^2}\left[-2\left(\frac{\phi'}{\phi}\right)'-n\left(\frac{\phi'}{\phi}\right)^2\right]+ \phi^{-2}R_{h}.
\end{equation}
\end{Lemma}

In our case, to compute the scalar curvature of $\bar g$, it suffices to substitute
$h = g_\Dcal + g_{\T^{n-3}}$, $\phi(\theta) = 1/z(r,\theta)$ and $a = \sinh r$ into \eqref{scalarWarpEq}. Noting that
$R_h = 1/2$, we have
\begin{align}
	R_{\bar g}&=(n-1)(\sinh r)^{-2}\bigl\{-2\partial_\theta[z\partial_\theta(1/z)]-n[z\partial_\theta(1/z)]^2\bigr\} + z^{2}/2\nonumber\\
	&=(n-1)(\sinh r)^{-2}\bigl\{2\bigl(\partial_\theta^2 z\bigr)/z-(n+2)[(\partial_\theta z)/z]^2\bigr\}	+ z^{2}/2,\label{Rgbar}
		\end{align}		
where $z$, $\partial_\theta z$ and $\partial_\theta^2 z$ are evaluated at $(r,\theta)$.

Now we are ready to observe the following.

\begin{Prop}\label{XrbarInfo}
For fixed $r>0$, the manifold $(\bar X_r, \bar g)$ satisfies:
\begin{enumerate}\itemsep=0pt
\item[$(a)$]  The scalar curvature
\[
	R_{\bar g}\geq \frac{1}{2}[z(r,\theta)]^2-C_{n,r}
\] for a
constant $C_{n,r}>0$ depending only on $n$ and $r$. In particular, we have
$
R_{\bar g}>-n(n-1)
$
provided that $p\in \H^2$ is chosen to have a large enough $z$-coordinate;

\item[$(b)$] Under the obvious identification $($isometry$)$ between $Y_{p,r}$ and $\bar Y_r$, we have
 $\bar H_r >H_{p,r}$ provided that the $z$-coordinate of $p$ is large enough.
 \end{enumerate}
\end{Prop}

\begin{proof}
	 $(a)$ follows from \eqref{Rgbar} and Lemma~\ref{dthetazEst};
	 $(b)$ follows from Lemma~\ref{YprMean}\,$(a)$ and \eqref{barHr}.
\end{proof}

\subsection{The gluing step}\label{gluingSec}

\begin{Lemma}[{\cite[Theorem 5]{BMN2011}}] \label{bmnthm5}
Let $\Omega$ be a compact $n$-manifold with boundary $\partial\Omega$, and let~$g$ and~$\tilde g$ be two smooth Riemannian metrics on $\Omega$ such that
\begin{enumerate}\itemsep=0pt
\item[$(a)$] $g-\tilde  g=0$ at each point on $\partial\Omega$;
\item[$(b)$] the mean curvatures satisfy $H_{\tilde g}-H_{g}>0$ at each point on $\partial \Omega$.
\end{enumerate}
Then, given any $\epsilon>0$ and any neighborhood $U$ of $\partial \Omega$, there exists a smooth metric $ \hat g$ on $\Omega$ with the following properties:
\begin{enumerate}\itemsep=0pt
    \item[$(1)$] $R_{\hat g} \geq \min\{R_g , R_{\tilde g} \}-\epsilon $ in $\Omega$;
	\item[$(2)$] $\hat g=\tilde g$ in $\Omega\setminus U$;
	\item[$(3)$] $\hat g= g$	in a neighborhood of $\partial\Omega$.
\end{enumerate}
\end{Lemma}
\begin{Remark}
By an arbitrary extension, in Lemma~\ref{bmnthm5} it suffices to assume that $g $ is defined only in a neighborhood of $\partial \Omega$. 	
\end{Remark}

To prove Theorem~\ref{counterexample},
a basic idea is to apply Lemma~\ref{bmnthm5} to obtain a
metric $\hat g$ on $\bar X_r$ which
agrees with $g_H$ in a neighborhood of $\partial \bar X_r\cong \partial X_{p,r}$,
so $\hat g$ extends smoothly into
$\bigl(\H^{n}/\Z^{n-2}\bigr)\setminus X_{p,r}$ by $g_H$.
A compromise is the $\epsilon$-cost to the scalar curvature estimate.
Thus, one would like to have a bit more scalar curvature
to begin with, so that the cost can be absorbed, maintaining the
desired lower bound $R_{\hat g}\ge -n(n-1)$.
This can be achieved by a suitable deformation of
$g_H$ in a neighborhood of $Y_{p,r}\subset \H^{n}/\Z^{n-2}$, as the following lemma
shows.

\begin{Lemma}\label{tradeoffLemma}
Let
\begin{equation}
u(\varrho)=
		\begin{cases}
			1-{\rm e}^{\frac{1}{\varrho-r_0}},& \varrho< r_0, \\
		    1, & \varrho\geq r_0,
		\end{cases}
\end{equation}
and define
\begin{equation*}%\label{metric2}
g_{H}':=\bigl[u(\varrho)\bigr]^2\ed\varrho^2 +\sinh^2 \varrho\, \ed\theta^2 +\bigl[z(\varrho,\theta)\bigr]^{-2} g_{\T^{n-2}}.
\end{equation*}
As long as $r_0>0$ is small enough, we can find $\delta>0$ such that
\[
	R_{g_{H}'} + n(n-1)  >0\qquad \mbox{for }
 \varrho \in [r_0-2\delta, r_0).
\]
\end{Lemma}

\begin{proof}

By \cite[Claim 2.1]{SWW2022}, we have
\begin{equation}\label{hatscalarcurv}
R_{g_H'}=R_{g_H}+\bigl(1-u^{-2}\bigr)(R_{\gamma(\varrho)}-R_{g_H})+2u^{-3}u'(\varrho) H_{p,\varrho},
\end{equation}
where
${\gamma(\varrho)}= \sinh^2 \varrho \,\ed\theta^2 +\bigl[z(\varrho,\theta)\bigr]^{-2}g_{\T^{n-2}}$ and $R_{g_H} = -n(n-1)$.

We want to estimate the right-hand side of formula~\eqref{hatscalarcurv}. To start with, by Lemma~\ref{scalarwarp},
\begin{equation*}%\label{Rgamma}
R_{\gamma(\varrho)}= (n-2) (\sinh \varrho )^{-2}
\big\{2\big(\partial_\theta^2 z\big)/z -(n+1)[(\partial_\theta z)/z]^2\big\}.
\end{equation*}
Thus, by the proof of Lemma~\ref{dthetazEst}, there exists a constant $C_{n,r_0}$,
depending on $n,r_0$ only, such that
\begin{equation}\label{RgammaBd}
|R_{\gamma(\varrho)}|\leq C_{n,r_0}\qquad\mbox{for } \varrho \in [r_0/2, r_0].
\end{equation}
Next, by the definition of $u$, we have, for $\varrho\le r_0$,
\begin{equation}\label{1minusu}
	0\ge 1 - u^{-2} = u^{-2}\bigl(-2 {\rm e}^{\frac{1}{\varrho - r_0}} + {\rm e}^{\frac{2}{\varrho - r_0}}\bigr) \ge - 2u^{-2} {\rm e}^{\frac{1}{\varrho - r_0}} \ge -2u^{-3} {\rm e}^{\frac{1}{\varrho - r_0}} .
\end{equation}
Moreover, for sufficiently small $r_0$, we have
$H_{p,\varrho} \geq 1$ for any $\varrho\le r_0$ (Lemma~\ref{YprMean}\,$(b)$), and so
\begin{equation}\label{urhoHprho}
	2u^{-3} u'(\varrho) H_{p,\varrho}  \ge 2u^{-3} {\rm e}^{\frac{1}{\varrho - r_0}}(\varrho - r_0)^{-2}.
\end{equation}
On combining \eqref{hatscalarcurv}, \eqref{RgammaBd}, \eqref{1minusu}
and \eqref{urhoHprho}, we obtain that
\begin{equation*}%\label{RgHprimeEst}
	R_{g_H'} - R_{g_H} \ge 2u^{-3} {\rm e}^{\frac{1}{\varrho - r_0}}
	\bigl[(r_0 - \varrho)^{-2} - C_{n,r_0} - n(n-1)\bigr]\qquad\mbox{for } \varrho\in [r_0/2, r_0].
\end{equation*}
Clearly, we can choose a small $\delta>0$ such that
\begin{equation*}%\label{RghDiffBD}
 R_{g_H'} - R_{g_H} >0 \qquad \mbox{for }
 \varrho \in [r_0-2\delta, r_0).
\end{equation*}
This completes the proof.
\end{proof}

\begin{figure}[t]
	\centering
	\includegraphics[scale = 0.5]{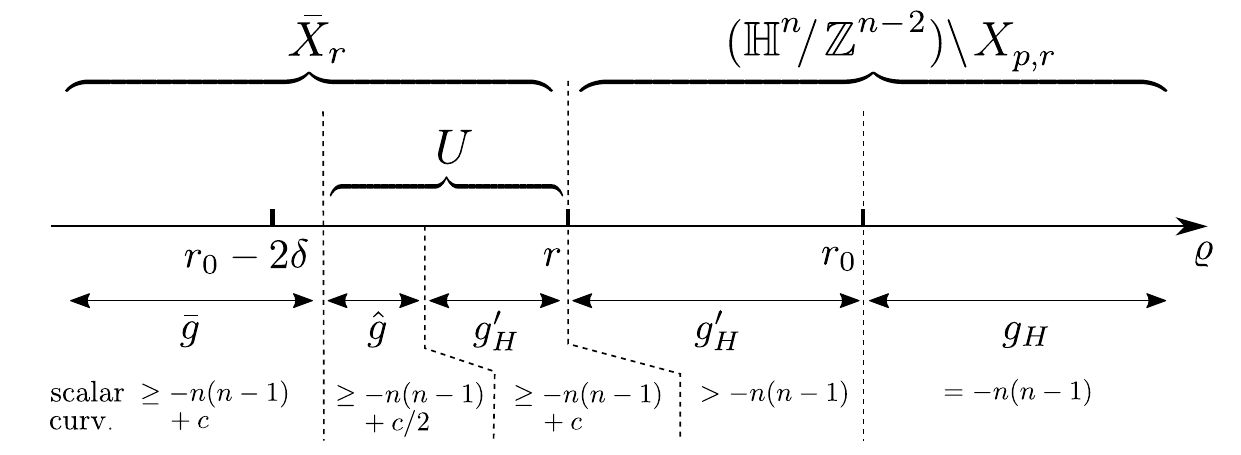}
	\caption{A schematic picture of $(\bs{M}, \bs{g})$.}\label{Fig_glue}
\end{figure}

\begin{proof}[Proof of Theorem \ref{counterexample}]
Let $r_0$ be small enough, and let $u(\varrho)$, $g_H'$ and $\delta$ be as in
Lemma~\ref{tradeoffLemma}. Define
\[
	c:= \min_{\varrho\in [r_0 - 2\delta, r_0 - \delta]} R_{g_H'} + n(n-1) >0.
\]	
Take $r := r_0 - \delta$, and note that we still have the freedom of
choosing $p\in \H^2$.

Suppose that the isometry between $Y_{p,r}$ and $\bar Y_r$
maps $\bs{q}\in Y_{p,r}$ to $\bar{\bs{q}}\in \bar Y_r$.
Furthermore, by using Fermi coordinates, any point in a small neighborhood of
$Y_{p,r}\subset X_{p,r}$ is uniquely represented by a pair
$(\bs{q}, d')$, where $d'$ is the $g_{H}'$-distance to $Y_{p,r}$.
Similarly, $(\bar{\bs{q}}, \bar d)$
coordinatizes a neighborhood of $\bar Y_r\subset \bar X_r$.
By identifying $(\bs{q}, d)$ with $(\bar {\bs q}, d)$, we have arranged that
$g_H' = \bar g$ along $\bar Y_{r}$.

To apply Lemma~\ref{bmnthm5}, assign $\Omega = \bar X_r$, $g = g_{H}'$ (defined
in a neighborhood $U$ of $\bar Y_r\subset\bar X_r$, via the identification above) and $\tilde g = \bar g$ (defined
on $\bar X_r$).
As noted above, Lemma~\ref{bmnthm5}\,$(a)$ is satisfied.
Furthermore, the mean curvature of $Y_{p,r}\subset X_{p,r}$ with respect to
$g_{H}'$ is $H_{p,r}':=H_{p,r}/u(r)\ge H_{p,r}$, but by choosing $p$ to have
large $z$-coordinate, we can still
arrange that $\bar H_r > H_{p,r}'$ (see the proof of Proposition~\ref{XrbarInfo}\,$(b)$).
Next, by shrinking $U$ if needed, we can assume that $R_{g_H'} \ge c - n(n-1)$ on~$U$, and we can assume the same lower bound for $R_{\bar g}$ by choosing $p$ suitably (Proposition~\ref{XrbarInfo}\,$(a)$).
Finally, take $\epsilon = c/2$.

With the above setting, Lemma~\ref{bmnthm5} applies and yields a metric
$\hat g$ defined on $\bar X_r$, satisfying
\begin{itemize}\itemsep=0pt
 	\item $R_{\hat g}\ge -n(n-1)+ c/2$;
	\item $\hat g = \bar g$ on $\bar X_r \setminus U$;
	\item $\hat g = g_H'$ in a neighborhood of $\bar Y_r\subset \bar X_r$.
\end{itemize}
Thus, $\hat g$ and $g_H'$ piece together to become a smooth metric
$\bs{g}$ defined  on
\[
	 \bs{M}:=\bigl[\big(\H^n/\Z^{n-2}\big)\setminus X_{p,r}\bigr] \cup \bar X_r/{\sim},
\]
where $\sim$
indicates boundary identification,
with (non-constant) scalar curvature
$R_{\bs{g}}\ge -n(n-1)$.
(For the reader's convenience, Figure~\ref{Fig_glue} includes a schematic,
1-dimensional illustration of the construction.)

In the statement of Theorem~\ref{counterexample},
take $(M,g) = (\bs{M}, \bs{g})$, $K = \bar X_r \cup (X_{p,r_0}\setdiff X_{p,r})\subset \bs{M}$
and $K' = X_{p,r_0}$, and the proof is complete.
\end{proof}

\subsection{Further remarks}
\subsubsection[Surgery applied to H\^{}n/Z\^{}(n-1)]{Surgery applied to $\boldsymbol{\H^n/\Z^{n-1}}$}\label{CErmksurg}
The construction above only modifies a portion of $\H^n/\Z^{n-2}$ that is contained in between $x_0< x < x_1$ for some
$x_0, x_1\in \R$. By a translation, we can always arrange that $x_0 = 0$.
Now, $T\colon (x,\bs{y},z)\mapsto (x+x_1, \bs{y},z)$ maps a neighborhood of $\{x = 0\}$ isometrically to a neighborhood of~$\{x = x_1\}$. Thus, by
removing the subsets $\{x<0\}$ and $\{x> x_1\}$ from $\bs{M}$ and then
identifying $\{x = 0\}$ and $\{x = x_1\}$  via $T$, we obtain a smooth Riemannian manifold $(\bs{M}', \bs{g}')$ that satisfies $R_{\bs{g'}}\ge -n(n-1)$. In fact, $(\bs{M}', \bs{g}')$ can be viewed as a
compactly supported `deformation' of~a~hyperbolic cusp $\H^n/\Z^{n-1}$, where $\Z^{n-1}$ acts on $(x,\bs{y})\in \R^{n-1}$ by translating
along the lattice $x_1\Z \times 2\pi\Z^{n-2}$. This serves as yet another counterexample to Gromov's Statement~\ref{generminoo}.

\subsubsection{A note on topology}\label{CErmktop}
It is interesting to determine the topology of both $\bs{M}$ and $\bs{M}'$ above.

Topologically, $\bs{M}$ is obtained by a \emph{surgery} along
$\{p\}\times \S^1 \subset \R^2\times \S^1$ and then taking product with $\T^{n-3}$.
The result of that surgery is homeomorphic to
$\S^1\times\R^2$. To see this, view $\S^3$
as the union of $\D^2\times \S^1$ and $\S^1\times \D^2$ with
the boundaries identified.
Then $\R^2\times \S^1$ is simply $\S^3$ with the core circle $\Ccal=\S^1\times \{q\}$ removed. Surgery of $\S^3$ along $\{p\}\times \S^1$ yields $\S^1\times\S^2$. Then removing~$\Ccal$ from $\S^1\times \S^2$ gives $\S^1\times \R^2$.
In conclusion, $\bs{M}\cong \S^1\times \R^2\times \T^{n-3}$, which is homeomorphic to $\H^{n}/\Z^{n-2}\cong \R^2\times \S^1\times \T^{n-3}$ via a map that switches the first two factors.

\begin{figure}[t]
	\centering
	\includegraphics[scale = 0.38]{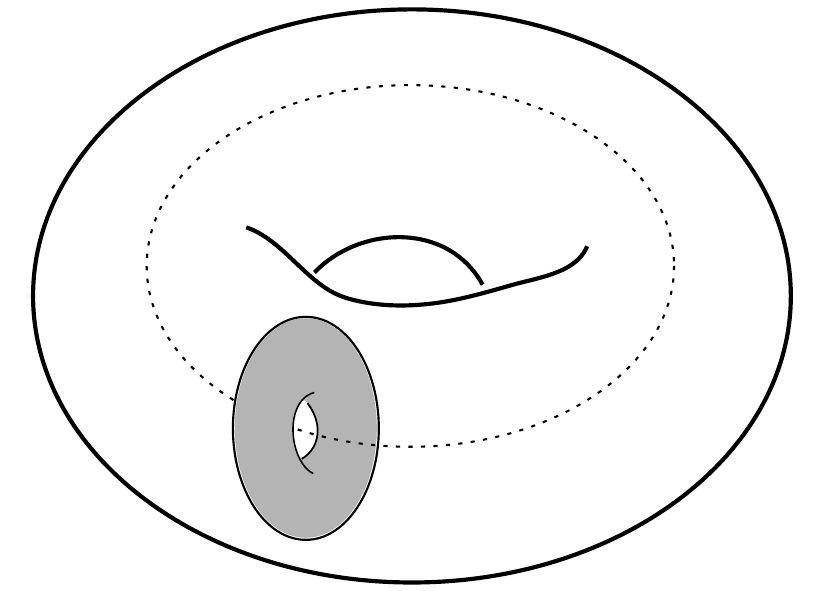}
	\caption{An illustration of $\big(\S^1\times \R^2\big)\setminus \big(\D^2_\epsilon\times \S^1\big)$,
	where $\D^2_\epsilon\times \S^1$ is shaded. }
	\label{Mprime}
\end{figure}

Regarding $\bs{M}'$, note that by identifying $x = 0$ and $x = x_1$
in $\{0\le x\le x_1\}\subset \H^2$, one obtains an open annulus, or equivalently $\R^2\setdiff\{\bs{0}\}$.
Thus, $\bs{M'}$ is obtained by a \emph{surgery} along $\{p\}\times \S^1\subset \big(\R^2\setdiff\{\bs{0}\}\big)\times\S^1$ and then taking product with $\T^{n-3}$.
In this case, a similar argument as the above applies, and the result of the surgery is homeomorphic to $\bigl(\S^1\times \R^2\bigr)\setminus \bigl(\D^2_\epsilon\times \S^1\bigr)$, i.e., the result of
removing a solid torus that is contained in a 3-ball $\B^3\subset \S^1\times \R^2$
(see Figure~\ref{Mprime}). Thus $\bs{M'} \cong \bigl[\bigl(\S^1\times \R^2\bigr)\setminus \bigl(\D^2_\epsilon\times \S^1\bigr)\bigr]\times \T^{n-3}$. In particular, the two ends of $\bs{M}'$ are separated
by a hypersurface with the topology $\S^2\times \T^{n-3}$; the same is \emph{not} true for $\H^n/\Z^{n-1}$.

\section{ALH manifolds, mass and deformations}\label{ALH_manifolds}
This section includes basic notions and results concerning
ALH manifolds, possibly with arbitrary ends, and their NRF and
conformal deformations. These results
 will be used in proving Theorem~\ref{alhpmt4}.

\subsection{ALH manifolds and mass}

\begin{Definition}\label{alh}
Let $(M^n,g)$ be a complete Riemannian manifold without boundary.
Suppose that
\begin{enumerate}\itemsep=0pt
    \item[$(1)$] for some (sufficiently large)
		compact set $K\subset M$, $M\setminus K$ has a connected component $\Ecal$
		that is diffeomorphic to $(0, 1)\times \mathbb{T}^{n-1}$, and
	\item[$(2)$] restricted to $\Ecal$, the metric $g$ admits an asymptotic expansion of the form
	\begin{equation}\label{alhasymp}
	g=\frac{1}{\tau^2}\left[\ed \tau^2+h+\frac{\tau^n}{n} \kappa +\Ocal\big(\tau^{n+1}\big)\right],
	\end{equation}
	where $\tau$ is the coordinate on the interval $(0,1)$; $h$ denotes a flat metric on $\mathbb{T}^{n-1} $, which represents metric at the conformal infinity $\Ecal_0$ (i.e., when $\tau = 0$); $\kappa=\kappa_{AB}\,\ed y^A \ed y^B$ is a symmetric tensor defined on $\T^{n-1}$, where $\bigl(y^A\bigr)$ are flat coordinates on $\T^{n-1}$; finally, $ \Ocal\bigl(\tau^{n+1}\bigr)$ stands for a remainder $Q = Q_{AB}\,\ed y^A\ed y^B$ with the asymptotics
	\begin{equation*}%\label{alhRmdr}
		\bigl|Q_{AB}\bigr| + \sum_{|\alpha|+k\le 2}\bigl|\tau^k\partial_y^\alpha \partial_\tau^k\  Q_{AB}\bigr| \le C\tau^{n+1} \qquad \mbox{as } \tau\rightarrow 0,
	\end{equation*}
	for some constant $C$, where $\alpha = (\alpha_1,\ldots, \alpha_{n-1})$ are multi-indices.
	\end{enumerate}
Such an $(M,g)$ is called \emph{asymptotically locally hyperbolic $($ALH\,$)$}, and $\Ecal$
 an \emph{ALH end}. Moreover, if $M\setminus \Ecal$ is non-compact,
we say that $(M,g)$ is \emph{ALH with arbitrary ends}.
\end{Definition}

\begin{Definition}[{cf.~\cite[Definition 1.1]{LN2015}}]\label{massDef}
Given a Riemannian manifold $(M,g)$ with an ALH end $\Ecal$ on which $g$ admits
the expansion \eqref{alhasymp}, we call
 \begin{equation*}%\label{massaspt}
 	m_{\Ecal,g}:=\tr_{h}\kappa =h^{AB} \kappa_{AB}
 \end{equation*}
the \emph{mass aspect function} associated to the pair $(\Ecal,g)$.
Furthermore, define
 \begin{equation*}%\label{massbar}
 \bar m_{\Ecal,g} :=\sup_{\mathbb{T}^{n-1}} m_{\Ecal,g}.
 \end{equation*}
 \end{Definition}

Throughout, let each $\tau$-level set in $\Ecal$ be denoted by $\Ecal_\tau$.
The following lemma shows how $\bar m_{\Ecal,g}$ is related to the mean curvature of $\Ecal_\tau\subset \Ecal$.

\begin{Lemma}\label{meancurvslice}
Let $\big(M^n,g\big)$ be a Riemannian manifold with an ALH end $\Ecal$. If $\bar m_{\Ecal,g}<0$, then there exist constants $\tau_0, C>0$ such that
\begin{equation}\label{mass_mean}
H_{\Ecal_\tau}\geq (n-1)+C\tau^n \qquad \mbox{for } \tau\le\tau_0,
\end{equation}
where $H_{\Ecal_\tau}$ is the mean curvature of $\Ecal_\tau$ computed
with respect to the `outward normal' $-\partial/\partial\tau$.
\end{Lemma}

\begin{proof}
Before making any assumption about $\bar m_{\Ecal,g}$, we have
\begin{equation}\label{meanCurvAsymp}
H_{\Ecal_\tau}=(n-1)-\frac{n-2}{2n}m_{\Ecal,g}\tau^n +\Ocal\big(\tau^{n+1}\big).
\end{equation}
For $\bar m_{\Ecal,g} < 0$, let us take
 $C = -\bar m_{\Ecal,g}/10$, and clearly \eqref{mass_mean} holds for some $\tau_0>0$.
\end{proof}

\subsection{NRF deformations}
Given a Riemannian $n$-manifold $(M^n, g_0)$, the \emph{normalized Ricci flow} (NRF), with initial metric~$g_0$, is by definition a smooth family of
Riemannian metrics $g(t)$ on $M$
satisfying the evolution equation
\begin{gather}\label{NRF}
\partial_t g =-2[\Ric_g+(n-1) g],\qquad
g(0) = g_0.
\end{gather}

  \begin{Lemma}\label{nrfdeform}
Suppose that $(M^n, g_0)$ is a complete Riemannian manifold with an ALH end $\Ecal$
that satisfies $R_{g_0}\ge -n(n-1)$ as well as the assumptions
$(\ref{curvbdAssu})$ and $(\ref{RbdAssu})$ in Theorem~{\rm \ref{alhpmt4}}.
Then there exists a $T>0$ such that, for $t\in (0,T]$, $g(t)$ is complete and satisfies \eqref{NRF} along with the following properties:
\begin{enumerate}[$(i)$]\itemsep=0pt
	\item $(\Ecal, g(t)|_{\Ecal})$ remains ALH, and $g(t)$ has the expansion $($see equation \eqref{alhasymp}$)$
			\begin{equation*}%\label{NRFalhexp}
				g(t) = \frac{1}{\tau^2}\left[\ed\tau^2 + h + \frac{\tau^n}{n}\kappa(t) + \Ocal\bigl(\tau^{n+1}\bigr)\right];
			\end{equation*}
	\item on $M$, $R_{g(t)}\geq -n(n-1)$ for all $t\in (0, T]$;
	\item if $g_0$ is not Einstein, then $R_{g(t)}> -n(n-1)$ for all $t\in (0,T]$; 	
	\item outside a compact subset in $M$, $R_{g(t)}\le -\alpha/2$ for $t\in (0,T]$;
	\item if $\kappa(0)=0$, then $\kappa(t)=0$ for all $t\in (0, T]$;
	\item if $\kappa(0)=0$, then for any $t\in (0,T]$ we have $R_{g(t)}+n(n-1)=\mathcal{O}(\tau^{n+1})$ as $\tau\rightarrow0$.
\end{enumerate}
\end{Lemma}

\begin{proof}
The existence of $g(t)$, $t\in (0,T]$, satisfying
\eqref{NRF} follows from the existence of a solution~$\tilde g(t)$, $t\in (0,\widetilde T]$, of the Ricci flow initiated at $g_0$.
They are related by a time-transformation:
\[
	g(t):={\rm e}^{-2(n-1)t}\tilde{g}\left(\Phi(t)\right),\qquad\mbox{where } \Phi(t) = \frac{{\rm e}^{2(n-1)t}-1}{2(n-1)}.
\]
Thus, up to constant factors, the curvature tensor $\Rm(t)$ of
$g(t)$ satisfies the same estimates as $\tilRm(\Phi(t))$ of $\tilde{g}(\Phi(t))$.
In particular, it follows from \cite{Shi1989} that, for all $t\in (0,T]$, $g(t)$ is complete, and $|\Rm(t)|$ is uniformly bounded.

Now we turn to proving the properties. $(i)$ follows from \cite[Proposition 3.1]{BW12}.
$(ii)$ can be verified by applying the maximum principle (see \cite[Theorem 7.42]{CLN06}) to
the evolution equation\footnote{For the evolution equation satisfied by $R_{g(t)}$, see \cite[formula~(5.1)]{BW12}.}
satisfied by ${\rm e}^{2(n-1)t}(R_{g(t)}+n(n-1))$;
to prove $(iii)$, invoke the strong maximum principle on the domain $\Omega\times[0,t]$,
 where $\Omega\subset M$ is compact on which $g_0$ is not Einstein,
 and then let $\Omega$ exhaust~$M$.
$(iv)$ would follow once we show that the integral
\begin{equation}\label{intRm}
	\int_{0}^t \partial_{t'}\tilRm\, \ed t', \qquad t\in \bigl(0,\widetilde T\bigr]
\end{equation}
is uniformly bounded;
to see this, note that the first covariant derivatives of $\tilRm(0)$ are assumed to be bounded (assumption $(d)$ in Theorem~\ref{alhpmt4}),
by \cite[Theorem~14.16]{CC2008}, we have
\[
\bigl|\nabla_{\tilde g(t)}^{2}{\tilRm}(t)\bigr|\leq\frac{C}{\sqrt{t}}
\]
for some constant $C>0$; in addition, the evolution equation of $\tilRm$ reads\footnote{${\tilRm}* {\tilRm}$ indicates a specific linear combination of the traces of $\tilRm\otimes\tilRm$.}
\[
\partial_{t}{\tilRm}=\Delta_{\tilde g(t)} {\tilRm}+{\tilRm}* {\tilRm};
\]
of course, $1/\sqrt{t}$ is integrable; combining these,
it is easy to see that \eqref{intRm} is uniformly bounded for small enough $\widetilde T$; since, by
assumption $(e)$ in Theorem~\ref{alhpmt4}, $R_g\le -\alpha$ outside a compact set, $(iv)$ follows.
$(v)$ follows from \cite[Proposition 4.3]{BW12}. Finally, $(vi)$ follows from $(v)$ and  \cite[formulas~(3.19)--(3.21)]{BW12}
(note that $g^{ij}(\tau)$ provides an extra factor of $\tau^2$).
\end{proof}

\subsection{Conformal deformations}
Throughout this section, $c_n:= 4(n-1)/(n-2)$.

\begin{Lemma}\label{confSubSol}
Let $(M,g)$ be complete with an ALH end $\Ecal$,
and let $f\in C^\infty(M)$ be a non-negative function that satisfies
\begin{enumerate}\itemsep=0pt
\item[$(a)$] $\supp f \subset K\cup \Ecal$ for some compact subset $K\subset M$;
\item[$(b)$] $f \in \Ocal(\tau^{n})$ as $\tau\rightarrow 0$ where $\tau$ is the function occurring in the expansion
		\eqref{alhasymp}.
\end{enumerate}
Then there exist a function $v\in C^\infty(M)$ and a constant $\delta_0$ such that  $0<\delta_0 \le v\le 1$ and
\begin{equation}\label{confLemma_veqn}
    -c_n \Delta_g v + fv = 0 \qquad\mbox{in } M.
\end{equation}
\end{Lemma}

\begin{proof}
    Let $\{\Omega_i\}_{i = 0}^\infty$ be a sequence of smooth, bounded domains satisfying $\Omega_i \Subset \Omega_{i+1}$ and $\bigcup_i \Omega_i = M$. For each $i$, the Dirichlet problem
    \begin{alignat}{3}
    &-c_n\Delta_g v_i + fv_i = 0\qquad &&\mbox{in } \Omega_i,&\nonumber\\
    &v_i= 1\qquad &&\mbox{on } \partial\Omega_i,&\label{viDirich}
    \end{alignat}
has a positive solution $v_i$. By the maximum principle, $0<v_i\le 1$. Thus, $v := \lim_{i\rightarrow\infty}v_i$ is well defined on $M$, satisfying $0\le v\le 1$ and \eqref{confLemma_veqn}. It remains to
show that $v$ has a positive lower bound.

Without loss of generality, assume that $\Sigma_i\subset \partial \Omega_i$ is the only component of $\partial \Omega_i$ that is contained in $\Ecal$; in fact, let us assume that each $\Sigma_i$ is a $\tau$-level set.
Denote $\tau_0:=\tau|_{\Sigma_0}$.

To refine the estimate of $v_i$, we construct an auxiliary function $\xi$ and compare it with $v_i$ via the maximum principle.
Indeed, let $\alpha\in (0,n-1)$ be any constant, and
define
\begin{equation*}%\label{auxXi}
    \xi =1 - (\tau/\tau_0)^\alpha, \qquad \tau\le \tau_0.
\end{equation*}
Using the fact that $-\ln \tau$ is, up to adding a constant, the distance
function to $\Sigma_0$, one easily computes that
\begin{equation}\label{LaplacianXi}
    \Delta_g \xi = \alpha (H_{\Ecal_\tau} - \alpha) ( \tau/\tau_0)^\alpha.
\end{equation}
Thus, by \eqref{meanCurvAsymp}, for sufficiently small $\tau_0$,
there exists a constant $C_{n,\alpha,\tau_0}>0$ such that
\[
	\Delta_g\xi \ge C_{n,\alpha,\tau_0} \tau^{\alpha}\qquad\mbox{for any } \tau\le \tau_0.
\]

Now, \eqref{viDirich}, the fact that $v_i\le 1$, and the assumption that
$f \in \Ocal(\tau^n)$ together imply
\begin{alignat*}{3}
            &\Delta_g v_i \le C'_{f,n} \tau^n\qquad &&\mbox{in } (\Omega_i\setdiff \Omega_0)\cap \Ecal,&\\
            &v_i > 0&&\mbox{on } \Sigma_0,&\\
            &v_i  = 1&&\mbox{on } \Sigma_i,&
\end{alignat*}
where $C'_{f,n}$ is a constant depending only on $f$ and $n$.
In comparison,
\begin{alignat*}{3}
            &\Delta_g\xi  \ge C_{n,\alpha, \tau_0} \tau^\alpha \qquad &&\mbox{in } \Ecal \setdiff \Omega_0,&\\
            &\xi = 0 && \mbox{on } \Sigma_0,&\\
            &\xi < 1 &&\mbox{on }  \Sigma_i.&
        \end{alignat*}
Thus, for sufficiently small $\tau_0$,
the maximum principle implies that $v_i\ge\xi$ in $(\Omega_i \setdiff \Omega_0)\cap \Ecal$.
Upon taking limit, $v\ge \xi >0$ on $\Ecal \setdiff \Omega_1$.
Since $v\ge 0$, the strong maximum principle, applied to \eqref{confLemma_veqn},
implies that $v>0$ on $M$.

When  $M\setminus \Ecal$ is compact (i.e., $M$ having no arbitrary end),
the above already implies that~$v$ has a positive lower bound.
When $M\setminus \Ecal$ is non-compact,
since $f$ is supported in $K\cup \Ecal$, by choosing $\Omega_0$ to include $K$, we have that each $v_i$ $(i\ge 1)$ is harmonic on $\Omega_i \setdiff(\Omega_0\cup \Ecal)$; using the maximum principle again, we get
\begin{equation*}
    \min_{\Omega_i \setdiff(\Omega_0\cup \Ecal)} v_i = \min_{\partial \Omega_0\setminus \Ecal} v_i
    \xrightarrow{i\rightarrow\infty} \min_{\partial\Omega_0\setminus\Ecal} v =: \delta_{\rm arb} >0.
\end{equation*}
To finish the proof, it suffices to take $\delta_0 = \min\{ \delta_{\rm arb} , \inf_{\Omega_1}v, \inf_{\Ecal\setminus\Omega_1}\xi\}$.
\end{proof}

\begin{Prop}\label{conformdeform1}
Let $(M^n, g)$ be complete, with an ALH end $\Ecal$ and with
$R_g \ge - n(n-1)$ on~$M$.
Let $\bar R\in C^\infty(M)$ be a function that satisfies
\begin{enumerate}\itemsep=0pt
\item[$(a)$] $-n(n-1)\leq \bar R\leq \min\{R_g, 0\}$;
\item[$(b)$] $\supp \bigl(R_g - \bar R\bigr) \subset \Ecal \cup K$ for some compact subset $K\subset M$;
\item[$(c)$]  $\bar R \equiv -n(n-1)$ on $\Ecal \setminus K'$ for some compact subset $K'\subset \Ecal$.
\end{enumerate}
Then the Yamabe equation
\begin{alignat}{3}
            &-c_n \Delta_g u  + R_g u  - \bar R u^{\frac{n+2}{n-2}} = 0\qquad&&\mbox{in } M,&\nonumber\\
            &u \rightarrow 1 \qquad&&\mbox{towards } \Ecal_0&\label{confyamabe}
        \end{alignat}	
has a solution $u$ with $ 0<\delta_0\le u \le 1$ for some constant $\delta_0$.
In particular, the metric $u^{4/(n-2)} g$ is complete and has the scalar
curvature $\bar R$.
\end{Prop}
\begin{proof}
The proof follows a super/sub-solution argument. To start with, define $L_g$ by
\begin{equation*}
    L_g u = -c_n \Delta_g u + R_g u - \bar R u^{\frac{n+2}{n-2}}.
\end{equation*}
Note that $L_g 1 = R_g - \bar R\ge 0$ by assumption. Thus, $1$ is a super-solution of \eqref{confyamabe}.

To find a sub-solution to \eqref{confyamabe}, take $f := R_g - \bar R\ge 0$.
Note that $R_g = - n(n-1)+ \Ocal(\tau^{n})$ in~$\Ecal$.
Thus, Lemma \ref{confSubSol} applies and yields a solution $v$ to \eqref{confLemma_veqn}, satisfying $0<\delta_0\le v\le 1$ for some constant $\delta_0$. Now we compute
\begin{equation*}
    -c_n\Delta_g v + R_g v - \bar R v^{\frac{n+2}{n-2}}
    = - c_n \Delta_g v + fv + \bar R\bigl(1- v^{\frac{4}{n-2}}\bigr)v = \bar R\bigl(1- v^{\frac{4}{n-2}}\bigr)v\le 0,
\end{equation*}
where the inequality follows from the assumption that $\bar R\le 0$ and the bounds for $v$. Thus, $v$ is a sub-solution of \eqref{confyamabe}.

Then one finishes the proof by following the argument of \cite[Proposition 2.1]{AM88}.
\end{proof}

Next, we will focus on the behavior of $u$
towards the ALH infinity $\Ecal_0$.

\begin{Lemma}\label{conformdeform2}
Let $u$ be as in Proposition $\ref{conformdeform1}$. Given any $\alpha\in (0, n-1)$, there exists a constant $\tau_0>0$ such that
\[
1- (\tau/\tau_0)^\alpha \leq u\leq 1 \qquad \mbox{for any }\tau\leq \tau_0.
\]	
\end{Lemma}

\begin{proof} Let $\xi:= 1 - (\tau/\tau_0)^\alpha$. By \eqref{LaplacianXi}, we have
\begin{equation}\label{LapXi2est}
	c_n\Delta_g \xi - \bigl(R_g - \bar R\bigr)\xi = c_n\alpha (H_{\Ecal_\tau} - \alpha) (\tau/\tau_0)^\alpha - \bigl(R_g - \bar R\bigr)\xi.
\end{equation}
Since $R_g-\bar R\in \Ocal(\tau^{n})$ in $\Ecal$, the right-hand side of \eqref{LapXi2est} is positive for $\tau\le \tau_0$, provided that $\tau_0$ is
sufficiently small.
On the other hand, since $\bar R \le 0$ and $0<u\le 1$, \eqref{confyamabe}
implies that
\begin{equation*}
	c_n \Delta_g u - \bigl(R_g - \bar R\bigr)u  \le 0.
\end{equation*}
Regarding boundary data,
\[
	u - \xi \ge 0 \qquad \mbox{along } \tau = \tau_0
	\qquad\mbox{and}\qquad	\lim_{\tau\rightarrow 0}u= \lim_{\tau\rightarrow 0} \xi = 1.
\]
Now the maximum principle implies that $u\ge \xi$ for $\tau \in (0,\tau_0]$.
\end{proof}

\begin{Prop}\label{behaviornearinfty}
Let $\big(M^n, g\big)$, $\bar R$ and $u$ be as in Proposition~{\rm \ref{conformdeform1}}.
Additionally, suppose that $R_g +n(n-1)\in \Ocal\bigl(\tau^{n+1}\bigr)$ and that, however small $\tau_0$ is, $R_g>-n(n-1)$ at some point in
 $\{\tau\leq \tau_0\}\subset \Ecal$. Then $u$ must have the following asymptotic expansion near $\tau = 0$:
\begin{equation*}%\label{expansionofu}
    u=1+u_{n0}\tau^n+\Ocal\bigl(\tau^{n+1-\epsilon}\bigr),
\end{equation*}
where $u_{n0}<0$ is a smooth function defined on the conformal infinity
$\Ecal_0\cong \T^{n-1}$ and $\epsilon >0$ is an arbitrary small constant.
\end{Prop}

\begin{proof} Let us take $w:=u- 1\leq 0$.
By \cite[Theorem 1.3]{ACF92}, $w$ has the expansion
 \begin{equation*}
     w=\sum\limits_{i=1}^{\infty}\sum\limits_{j=0}^{N_{i}}u_{ij}\tau^i(\log\tau)^j,
 \end{equation*}
where $u_{ij}\in C^{\infty}(\Ecal_0)$.
Clearly, the proof would be complete once we verify the conditions:
\begin{enumerate}\itemsep=0pt
    \item[$(C1)$] $u_{ij} =0$ for $i<n$;
    \item[$(C2)$] $u_{nj} =0$ for  $j>0$;
    \item[$(C3)$] $u_{n0} <0$.
\end{enumerate}

\emph{Verification of $(C1)$.}
By \eqref{confyamabe}, $w$ satisfies
\begin{equation*}
	\Delta_g w  - n w = \frac{1}{c_n} \left[R_g (w+1)  - \bar R (w+1)^{\frac{n+2}{n-2}}\right]
						- nw.
\end{equation*}
Since only a neighborhood of $\Ecal_0$ is
concerned, we can simply substitute $\bar R = -n(n-1)$; by rearranging terms,
we get
\begin{equation*}%\label{LapminusnW}
	\Delta_g w - nw = \frac{1}{c_n} [R_g+n(n-1)] u
			 + \frac{n(n-1)}{c_n} \left[(w+1)^{\frac{n+2}{n-2}} - 1 - \frac{n+2}{n-2}w\right] =: A+B.
\end{equation*}
Since $\lim_{\tau\rightarrow 0 } u = 1$ and
$0\le R_g+n(n-1) \in \Ocal\big(\tau^{n+1}\big)$, we have $A\ge 0$ and $A\in \Ocal\big(\tau^{n+1}\big)$.
On the other hand,
$B$ is the remainder of a Taylor expansion truncated at the linear term, so
$B = \Ocal(w^2)$ as $\tau\rightarrow 0$. By Lemma~\ref{conformdeform2}, $w = \Ocal(\tau^\alpha)$ for any $\alpha< n- 1$.
Of course, we can choose $\alpha > (n+1)/2$, and thus
$B = \Ocal\bigl(\tau^{2\alpha}\bigr) = o\bigl(\tau^{n+1}\bigr)$.
In summary, for sufficiently small $\tau_0$, we have
\begin{equation}\label{LapwminusNwEst}
	0\le  \Delta_g w - n w \in \Ocal\bigl(\tau^{n+1}\bigr)\qquad \mbox{for } \tau \le \tau_0.
\end{equation}

Now consider any $\beta \in (n-1, n)$. Using \eqref{meanCurvAsymp},
it is easy to verify that
\[
	\Delta_g \tau^\beta - n\tau^\beta = -(\beta + 1)(n - \beta) \tau^\beta
	+ \Ocal\bigl(\tau^{n+2}\bigr).
\]
Clearly, there exists $\tau_0>0$ such that
\begin{equation*}
(\Delta_g - n)\big(w + \lambda \tau^\beta\big) \le 0 \qquad \text{for all }
	\tau\le \tau_0 \text{ and constants } \lambda \ge 1.
\end{equation*}
Fix such a $\tau_0$, and let us choose $\lambda\ge 1$ such that
$w|_{\{\tau = \tau_0\}}+\lambda\tau_0^\beta \ge 0$; moreover,
we have $\lim_{\tau\rightarrow 0} (w+\lambda\tau^\beta) = 0$.
Thus, by the maximum principle,
\begin{equation*}
	w\ge -\lambda\tau^\beta\qquad \mbox{for }\tau\le \tau_0.
\end{equation*}
Since $\beta\in (n-1,n)$ is arbitrary and $w\le 0$, this verifies $(C1)$.

\emph{Verification of $(C2)$.}
By $(C1)$, we have
\[
	w = \sum_{j = 0}^{N_n} u_{nj} \tau^n(\log\tau)^j + \Ocal\bigl(\tau^{n+1 - \epsilon}\bigr).
\]
Further information about $u_{nj}$ is obtainable by
computing $(\Delta_g - n)w$ using this expansion and then comparing the result with \eqref{LapwminusNwEst}. In fact, direct computation and \eqref{meanCurvAsymp}
yield:
\begin{gather*}
	(\Delta_g  - n)\tau^n  = \Ocal\bigl(\tau^{n+2}\bigr), \\
	(\Delta_g - n) [\tau^n(\log\tau)^j]
	 = \left[(n+1)j (\log\tau)^{j-1} + j(j-1)(\log\tau)^{j-2}\right]\tau^n+\Ocal\bigl(\tau^{n+2}\bigr)
\end{gather*}
with $1\le j\le N_i$.
Now, since $u_{nj}$ are all defined on $\Ecal_0$,
we have $\Delta_g u_{nj}\in \Ocal\big(\tau^2\big)$;
and since the remainder $\Ocal(\tau^{n+1-\epsilon})$ does not contribute to the coefficients $s_j$ of $\tau^n(\log\tau)^j$ in $(\Delta_g - n)w$, we have
that $s_j$ equals to
\begin{alignat*}{3}
		&(n+1)  u_{n1} + 2 u_{n2}\qquad &&\mbox{for } j = 0,&\\
		&2(n+1) u_{n2} + 6 u_{n3} \qquad &&\mbox{for } j= 1,&\\
		&\vdots&&&\\
		&N_n (n+1) u_{nN_n} \qquad &&\mbox{for }j = N_n - 1,&\\
		&0 \qquad &&\mbox{for }j = N_n.&
	\end{alignat*}
By \eqref{LapwminusNwEst}, all $s_j$ must vanish, which implies that
\begin{equation*}
	u_{nj} \equiv 0 \qquad\mbox{for } j =1, \ldots, N_n.
\end{equation*}
This verifies $(C2)$.

\emph{Verification of $(C3)$.}
Consider an auxiliary function $\zeta := -\delta \big(\tau^n + \tau^{n+1}\big)$
where $\delta>0$ remains to be chosen.
Now
\[
	(\Delta_g - n) \zeta = - \delta\big[(n+2) \tau^{n+1} + \Ocal\big(\tau^{n+2}\big)\big],
\]
so $(\Delta_g - n)\zeta \le 0$ provided that $\tau$ is small, and let us choose
$\tau_0$ accordingly (note: this is
independent of the choice of $\delta$).
By comparison, recall from \eqref{LapwminusNwEst} that
$(\Delta_g - n) w\ge 0$ for~$\tau\le \tau_0$.

Regarding boundary data, first note that the assumption about $R_g$
implies that $w$ cannot be identically zero for $\tau\in (0,\tau_0]$; thus, the strong maximum principle
implies, in particular, that $w < 0$ along $\tau = \tau_0$. This allows us to choose
$\delta$ such that $w \le \zeta$ along $\tau = \tau_0$.
Moreover, both $w, \zeta\rightarrow 0$ as $\tau\rightarrow 0$. Now, by the maximum
principle, we get
\[
	w\le \zeta = -\delta\bigl(\tau^n + \tau^{n+1}\bigr) \qquad\mbox{for } \tau\le \tau_0.
\]
This proves that $u_{n0}<0$, verifying $(C3)$.
\end{proof}

\begin{Lemma}\label{confdeformmass}	
	Let $(M^n,g)$ be a Riemannian manifold with an ALH end $\Ecal$,
	on which the asymptotic expansion \eqref{alhasymp} applies.
	Suppose that $u = 1+\varphi \tau^n + \Ocal(\tau^{n+1})$
	is a function defined on $\Ecal$, where $\varphi\in C^\infty(\Ecal_0)$.
	Then, up to a diffeomorphism that restricts to be the identity on $\Ecal_0$, the
	deformed metric $\bar g = u^{\frac{4}{n-2}} g$ on $\Ecal$
	has the expansion
	\begin{equation*}%\label{bargExp}
		\bar g = \frac{1}{\bar\tau^2}
			\left[ \ed\bar\tau^2 + \bar h + \frac{\bar\tau^n}{n}\bar\kappa
			+ \Ocal\bigl(\bar\tau^{n+1}\bigr)\right],
	\end{equation*}
	where
	\begin{equation*}
		\bar h = h\qquad\mbox{and}\qquad
		\bar \kappa = \kappa + \frac{4(n+1)}{n-2}\varphi h.
	\end{equation*}	
\end{Lemma}
\begin{proof}
	A standard argument following the proof of  \cite[Lemma 6.5]{BQ08}.
\end{proof}

\section{Two rigidity results}\label{Two_rigibility}
The goal of this section is to prove Theorem~\ref{alhpmt4}, Corollary~\ref{hpbInfCor} and Theorem~\ref{deform1}.
The reader may consult Appendix~\ref{mubbSec} before proceeding.

\begin{Prop}[{cf. \cite[Theorem~1.1]{CLSZ2021}}]\label{incprssArg}
    For $3\leq n\leq 7$, let $M^n$ be a $($connected$)$ non-compact manifold with connected, compact boundary $\Sigma$.
    Let $\iota\colon \Sigma\hookrightarrow M$ be the inclusion map. Suppose that $\Sigma\in \mathcal{C}_{\deg}$ $($see Definition~{\rm \ref{cdegdef})}
    and that the map
    $\iota$ is incompressible. Then $M$ admits no complete metric $g$ with  $R_g\geq-n(n-1)$ and $H_{\Sigma}> n-1$.
\end{Prop}
\begin{proof}
    To begin with, by the classification of covering spaces, there exists a covering of $M$, say $p\colon\hat{M}\rightarrow M$, that satisfies
    \begin{equation} \label{fundgrouopRel}
    	p_*\bigl(\pi_1\bigl(\hat{M}\bigr)\bigr)=\iota_*(\pi_1(\Sigma))\subset \pi_1(M),
    \end{equation}
    where base points for the fundamental groups are omitted. Moreover, by the homotopy lifting property, there exists an embedding
    $\hat{\iota}\colon \Sigma \rightarrow \hat{M}$ such that $\iota = p\circ\hat \iota$.

    By \eqref{fundgrouopRel} and the incompressibility of $\iota$, the composition
    \[
    	J:= \big({\iota_*}^{-1}\big|_{\iota_*(\pi_1(\Sigma))}\big)\circ p_*\colon \  \pi_1(\hat{M})\rightarrow \pi_1(\Sigma)
    \]
    is a well-defined group homomorphism.
     Since $\Sigma$ is aspherical, by \cite[Proposition 1B.9]{Hatcher_AT}, there exists a map $j\colon \hat{M}\rightarrow \Sigma$ such that $j_*\colon \pi_1(\hat{M})\rightarrow \pi_1(\Sigma)$ is equal to $J$;
     in particular,  $j_*\circ\hat\iota_* = \id_{\pi_1(\Sigma)}$; then, by applying the uniqueness part of  \cite[Proposition 1B.9]{Hatcher_AT} to $\Sigma$, it is easy to see that
     $j\circ \hat\iota$ is in fact homotopic to $\id_{\Sigma}$.

	Since $\iota$ is an embedding, each boundary component of $\hat M$, which is a lifting of $\Sigma$, must be
	diffeomorphic to $\Sigma$. In particular, denote $\hat{\Sigma}=\hat\iota(\Sigma)$.
	Since $j\circ \hat \iota$ is homotopic to $\id_\Sigma$, we have $\bigl[\hat{\Sigma}\bigr] = \hat\iota_*[\Sigma]\neq 0 \in H_{n-1}(\hat{M};\mathbb{Z})$.
		
	Now, for the sake of deriving a contradiction, suppose that $g$ is a complete metric on $M$ with $R_g\geq-n(n-1)$ and $H_{\Sigma}\geq (n-1)(1+\delta)$ for some constant $\delta>0$. Let $\hat{g}:=p^*g$ be the pull-back metric on $\hat M$, and define
	$\rho(x):=\dist_{\hat{g}}(x,\hat{\Sigma})$ for $x\in \hat{M}$.
	
	For an arbitrarily large $T>0$,  let
	\[\mathcal{D}_T:=\bigl\{x\in \hat{M}\colon \rho(x)\leq T\bigr\},\]
	and let $\hat{\Sigma}_i$ $(0\leq i\leq k)$ be those components of $\partial{\hat M}$ that satisfy
	\[
	\hat{\Sigma}_i \cap \mathcal{D}_T\neq \varnothing,
	\]
	where $\hat \Sigma_0 = \hat\Sigma$.
	Define  (see Figure~\ref{UcalTeps} below)
	\[
	\mathcal{U}_T= \mathcal{D}_T\cup \bigcup_{0\leq i\leq k} \hat{\Sigma}_i\qquad \text{and} \qquad
	\mathcal{U}_{T,\epsilon}=\bigl\{x\in\hat{M}\colon \dist_{\hat{g}}(x,\Ucal_T)<\epsilon\bigr\}.
	\]
	Since $M$ is complete, connected and non-compact, so is $\hat M$, and we have $\bar\Ucal_T \Subset \Ucal_{T,\epsilon}$.
	Moreover, for small enough $\epsilon$,
	\[
	\mathcal{U}_{T,\epsilon} \cap \bigg(\partial\hat{M}- \bigcup_{0\leq i\leq k}\hat{\Sigma}_i\bigg)=\varnothing.
	\]
	Thus, by the smooth Urysohn lemma, there exists a function $\eta \in C^\infty(\hat M)$ with
	\begin{equation*}
		\eta(x)=
		\begin{cases}
0,&  x\in \mathcal{U}_T,\\
1,&   x\in \hat{M}\setminus\mathcal{U}_{T,\epsilon}.
		\end{cases}
	\end{equation*}
	Let $a\in (0,1)$ be a regular value of $\eta$. Automatically, $\eta^{-1}(a)$ is a smooth, closed hypersurface of~$\hat M$, and
	$\eta^{-1}(a)\cap \partial\hat M = \varnothing$.

            \begin{figure}[t]
            \centering
            \includegraphics[scale = 0.25]{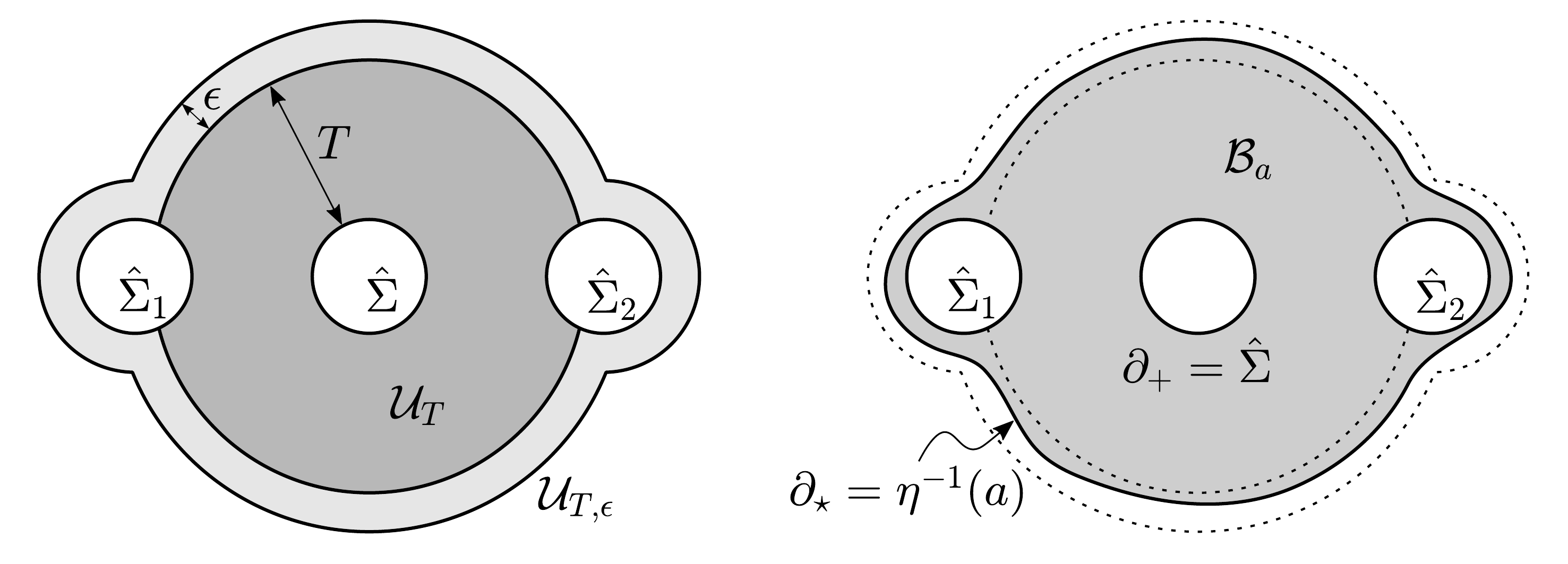}
            \caption{A schematic picture showing
            $\Ucal_T$, $\Ucal_{T,\epsilon}$ (left figure) and
            $\Bcal_a$ (right figure).
            The complement of~$\Ucal_{T,\epsilon}$, which may include more boundary components of~$\hat M$,
            is not displayed.}
            \label{UcalTeps}
            \end{figure}

	By the above arrangement, $\Bcal_a:=\eta^{-1}\bigl([0,a]\bigr)$, equipped with the restriction of the metric $\hat g$, is a Riemannian
	band with
	\[
	\partial_{+}=\hat{\Sigma}\qquad \mbox{and} \qquad
	\partial_{-}=\partial \Bcal_a \setminus \hat{\Sigma} = \eta^{-1}(a) \cup \bigcup_{1\le i\le k} \hat\Sigma_i .
	\]
	By letting $f = j|_{\Bcal_a}$ and using Lemma~\ref{C_degNSepLemma}, one easily sees that $\Bcal_a$ satisfies the \Nsep property (see Definition~\ref{noSepDef}).
	Then take $\partial_{\star}=\eta^{-1}(a)$.
		
	With these choices, all assumptions of Lemma \ref{bandWidth} are satisfied for $(\Bcal_a, \hat g|_{\Bcal_a};\partial_-, \partial_+)$
	and $\partial_\star$.
	Since $\bigl(\hat M,\hat g\bigr)$ is complete and non-compact,
	the distance $\dist_{\hat{g}}(\partial_{\star},\partial_{+})$ can get arbitrarily large as one chooses large $T$.
	This contradicts Lemma \ref{bandWidth}.
\end{proof}

\begin{Remark}
Proposition \ref{incprssArg} still holds if $M$ is allowed to be compact.
In fact, proceeding along the same proof, we still have \smash{$\big[\hat \Sigma\big]\ne 0\in H_{n-1}\big(\hat M;\Z\big)$}, so
$\hat M$ cannot be compact with a~single boundary component.
Hence, either (1) $\hat M$ is non-compact, and the previous proof applies verbatim;
or (2) $\hat M$ is itself a Riemannian band with $\partial_+ = \hat\Sigma$ that
satisfies the \Nsep property and the curvature bounds $R_{\hat g}\ge -n(n-1)$,
$H_{\partial \hat M}\ge (n-1)(1+\delta)$; however, by
Remark~\ref{bandNonexistence}\,$(A)$, such a band cannot exist, reaching a contradiction.
\end{Remark}

\begin{Prop}\label{alhPrep}
	For $3\le n\le 7$, let $\big(M^n,g\big)$ be a complete Riemannian manifold without boundary, with an ALH end
	$\Ecal\cong (0,1)\times \T^{n-1}$, and satisfying $R_g\ge -n(n-1)$.
	Suppose that $Y:=M\setminus \Ecal$ is non-compact and that $\partial Y\cong \T^{n-1}$
	is incompressible in $M$. Then $\bar m_{\Ecal, g}\ge 0$.
	In~addition, if the assumptions $(\ref{curvbdAssu})$, $(\ref{RbdAssu})$ in Theorem~{\rm \ref{alhpmt4}} hold,
	then $\kappa = 0$ only if $(M,g)$ is Einstein.
\end{Prop}

\begin{proof}
Suppose, on the contrary, that $\bar m_{\Ecal,g}< 0$.
Let $\tau$ be a defining function compatible with the ALH structure of $\Ecal$
(see \eqref{alhasymp}).
Then by Lemma \ref{meancurvslice},
there exists a small $\tau_0>0$ such that the mean curvature of the level set $\Ecal_{\tau_0}$ satisfies
\[
	H_{\Ecal_{\tau_0}} \geq (n-1)+\delta_0
\]
for some $\delta_0>0$.

Now, remove
$\{0< \tau< \tau_0\}$, a subset of $\Ecal$, from $M$ and denote the resulting manifold by $M'$.
By using the assumptions, it is easy to see that $\partial M' = \Ecal_{\tau_0}\cong \T^{n-1}$ is incompressible in $M'$.
Clearly, $\partial M'\in \Ccal_{\deg}$.
By Proposition~\ref{incprssArg},
we get a contradiction.
This proves the inequality $\bar m_{\Ecal, g}\ge 0$.

Next we turn to the second part of the proposition. Again we argue by contradiction.
Assume that $\kappa = 0$
without $(M,g)$ being Einstein. Let $g(t)$ be the NRF initiated at $g$.
Then by Lemma~\ref{nrfdeform}, for some small $t_0$, we have
\begin{enumerate}[$(i)$]\itemsep=0pt
	\item $R_{g(t_0)}>-n(n-1)$ on $M$;\label{scalLB}
	\item $R_{g(t_0)}\le -\alpha/2<0$ outside a compact subset of $M$;
	\item $R_{g(t_0)} = -n(n-1) + \Ocal\big(\tau^{n+1}\big)$ on $\Ecal$;	\label{scalAsympE}
	\item $\bigl(\Ecal, g(t_0)|_\Ecal\bigr)$ remains ALH with $\kappa(t_0) = 0$.
\end{enumerate}
It is easy to check that a function $\bar R$ as described in Proposition~\ref{conformdeform1} exists; thus, there is
a positive function $u\in C^\infty(M)$ such that $\bar g := u^{4/(n-2)}g(t_0)$ is complete
with  $R_{\bar g} = \bar R\ge -n(n-1)$. Furthermore, thanks to ($\ref{scalLB}$) and ($\ref{scalAsympE}$) above,
both Proposition~\ref{behaviornearinfty} and Lemma~\ref{confdeformmass} apply.
As a~consequence, $\bigl(\Ecal, \bar g|_\Ecal\bigr)$ remains ALH and satisfies
\[ \bar\kappa =  \frac{4(n+1)}{(n-2)} u_{n0} h,\] where
$u_{n0}<0$, and $h$ is a flat metric on $\T^{n-1}$.
Clearly, $\bar m_{\Ecal, \bar g} =\tr_h\bar\kappa < 0$. This contradicts the first part of
 the proposition.
\end{proof}

\begin{proof}[Proof of  Theorem \ref{alhpmt4}]	
	For convenience, let $\Ncal_o$ (resp., $\Hcal_o$) denote the
	result of removing a~tubular neighborhood of $\phi\big(\T^k\big)$ from $N$ (resp., $\psi\big(\T^k\big)$ from $\H^n/\Z^{n-1}$). Both
	$\partial \Ncal_o$ and $\partial \Hcal_o$ inherit the product structure $\S^{n-k-1}\times \T^{k}$, which are
	identified to form $M$. In symbols, $M = \Hcal_o\sqcup_\Phi \Ncal_o$, where $\Phi\colon  \partial \Hcal_o\rightarrow \partial \Ncal_o$
	is the identification map.
		
	By Proposition~\ref{alhPrep}, to prove the theorem it suffices to show that the boundary
	$\Sigma$ of $M\setminus\Ecal$ is incompressible in $M$.
	
	To show this, it in turn suffices to show that the $\T^k$-factor of $\partial \Ncal_o$
	is incompressible in $M$, according to Lemma~\ref{incprLemma}.
	
	If this was not the case, let $L$ be a non-contractible loop in
	$\{x\}\times \T^{k}\subset \partial \Ncal_o$
	that is contractible in $M$.
	
	Now consider
	\[	
		\Hcal':=\bigl(\S^1\times \T^{n-k-1} - \Bb\bigr)\times \T^k,
	\]
	where $\Bb$ is an $(n-k)$-ball embedded in $\S^1\times \T^{n-k-1}$.
	Topologically, $M$ can be viewed as a~subset of $M':=\Hcal' \sqcup_\Phi \Ncal_o$, so $L$ is also contractible
	in $M'$. By \cite[Lemma~{A.3}]{CLSZ2021}, $\Hcal'$ satisfies the `lifting property' (see \cite[Definition~{A.2}]{CLSZ2021}).
	Thus, \cite[Lemma~{A.4}]{CLSZ2021} applies, showing that $L$ is contractible in $\Ncal_o$ and hence in $N$;
	since $\{x\}\times\T^{k}$ and $\phi(\T^k)$ are homotopic in $N$, $\phi$ cannot be incompressible, violating the assumption
	$(\ref{incprsAssu})$.
\end{proof}
\begin{Remark}\label{vanKampenRmk}
The proof above can be made more direct if one assumes that $k<n-2$. In this case, both $\pi_1(\partial \Ncal_o)$ and $\pi_1(\partial \Hcal_o)$ are isomorphic to $\pi_{1}\big( \T^{k}\big)$, and it is easy to see that the maps $\pi_1\left(\partial \mathcal{N}_o\right) \rightarrow \pi_1\left(\mathcal{N}_o\right)$ and $\pi_1\left(\partial \mathcal{H}_o\right) \rightarrow \pi_1\left(\mathcal{H}_o\right)$
are both injective.
By van Kampen's theorem, we have $\pi_1(M) \cong \pi_1\left(\mathcal{H}_o\right) *_{\pi_1\left(\partial \mathcal{N}_o\right)} \pi_1\left(\mathcal{N}_o\right)$. Thus, a direct
application of
\cite[Theorem~11.67\,$(i)$]{Rotman_group} shows that
$\partial\mathcal{N}_o$ is incompressible in $M$, and it follows that the $\T^k$-factor of $\partial\mathcal{N}_o$ is also incompressible in $M$.
\end{Remark}

\begin{proof}[Proof of Corollary~\ref{hpbInfCor}] In this setting, the assumptions $(a-e)$
in Theorem~\ref{alhpmt4} are satisfied. Since $\kappa$ automatically vanishes, we conclude that $g$ is Einstein.
Write the metric on $\H^n/\Z^{n-1}$ as $\ed t^2 + {\rm e}^{2t} g_0$ where $g_0$ is a flat metric on $\T^{n-1}$. Since $\H^n/\Z^{n-1}$ is isometric
to $(M,g)$ outside a compact set, one can remove the corresponding cusp (i.e., $\{t<-a\}$ for some $a\gg 0$) from $M$ and obtain a
complete, non-compact manifold $(M', g')$ with boundary $\partial M'\cong \T^{n-1}$, satisfying $H_{\partial M'} \equiv  n - 1$, where the mean
curvature is computed with respect to the \emph{inward} normal.
By \cite[Theorem 2]{CK92}, $(M', g')$ is isometric to $[-a,\infty)\times \T^{n-1}$ with the warped product metric
$\ed {t'}^2 + {\rm e}^{2t'} g_0$;
by using this fact and the respective distance functions to $\partial M'\subset M$ and $\{-a\}\times\T^{n-1}\subset \H^n/\Z^{n-1}$, it is easy to construct an isometry
between $(M,g)$ and $\H^n/\Z^{n-1}$.
\end{proof}

\begin{Remark}\label{hpbInfCor_NonCptRmk}
The statement of Corollary~\ref{hpbInfCor} remains true when $N$ is
non-compact without boundary. In fact, one only needs to prove the incompressibility of
a $\T^{n-1}$-slice located far into the ALH infinity of $M$, and this is handled by a corresponding step
in the proof of Theorem~\ref{alhpmt4}. Then the result follows directly from Theorem~\ref{alhpmt1}.
\end{Remark}

\begin{proof}[Proof of Theorem~\ref{deform1}]
Let $(x,z)$ be the standard coordinates on $\R^2_+$, a topological factor of~$\H^n/\Z^{n-2}$. Since $K$ is compact, via the isometry $f$,
both $x$ and $z$ can be regarded as coordinate functions on $M\setminus K$. Thus, for a large enough $x_0>0$, we can
remove $\{|x|>x_0\}$ from $M$ and then identify $\{x = \pm x_0\}$ in the same way
as we did in Section~\ref{CErmksurg}. The result is a~complete Riemannian manifold $(M^*,g^*)$
with an ALH end $\Ecal$, satisfying $R_{g^*}\ge -n(n-1)$. Moreover, $(M^*, g^*)$ is isometric to $\H^n/\Z^{n-1}$ outside a compact set; thus,
the assumptions ($\ref{curvbdAssu}$), ($\ref{RbdAssu}$) in Theorem~\ref{alhpmt4} hold automatically, and
$\kappa = 0$ for $(\Ecal, g^*|_\Ecal)$.

It is easy to see that $M^*$ is of the form $M_1\sqcup_\Phi M_2$ as
described in Lemma~\ref{incprLemma}
with $k = n-2$. In particular, $M_2$ can be viewed as a subset of $M$.
By assumption, $f^{-1}(T)$ is incompressible in~$M$ and hence in~$M_2$.
Using the proof of Theorem~\ref{alhpmt4}, one can show that $f^{-1}(T)$ is incompressible in $M^*$;
then by Lemma~\ref{incprLemma}, $\partial (M\setminus\Ecal)\cong\T^{n-1}$ is incompressible
in $M^*$.

Thus, all conditions in Proposition~\ref{alhPrep} are verified for $(M^*, g^*)$, and we conclude that $g^*$ is Einstein.
The proof of Corollary~\ref{hpbInfCor} shows that there is an isometry
$\tilde f\colon (M^*, g^*)\rightarrow \H^n/\Z^{n-1}$ that uniquely extends the isometry, induced by $f$, between the `cuspidal ends'
%$\{0<z<z_0\}$
in $M^*$ and $\H^n/\Z^{n-1}$.
Let $z_0>0$ be sufficiently small; then by using distance functions to the hypersurfaces
$\{z = z_0\}$ in both $M$ and $\H^n/\Z^{n-2}$, it is easy to
construct an isometry between $(M,g)$ and $\H^n/\Z^{n-2}$; details are
left to the interested reader.
\end{proof}

\section{Two splitting results of `cuspidal-boundary' type}\label{Two_splitting}
The bulk of this section is dedicated to proving Theorem~\ref{alhpmt2}.
The proof of
Theorem~\ref{alhpmt1}, which largely depends on those of Proposition~\ref{incprssArg}
%(for the inequality part)
and Theorem~\ref{alhpmt2},
%(for the `splitting' part)
will be sketched at the end
of the section.

Now we begin our proof of Theorem~\ref{alhpmt2}.

In addition to its hypothesis, let us assume that $H_{\partial M} \ge 3$.
Under this assumption, the proof would be complete once we show that
$(M,g)$ is isometric to $\bigl((-\infty,0]\times \Sigma, \ed t^2 + {\rm e}^{2t} g_0\bigr)$ for some
closed $3$-manifold $\Sigma$ carrying a flat metric $g_0$. In fact,
 $\Sigma$ will occur as a hypersurface in $M$, obtained by an approximation scheme involving $\mu$-bubbles (Sections~\ref{tauketc} and \ref{muk_in_Ek}); then we show that $\Sigma$ must be compact and that $(M,g)$ is isometric to the desired warped product (Section~\ref{Sigmak_limits}).

The reader is recommended to consult Appendix~\ref{mubbSec} before proceeding.

\subsection[Specification of mu\_{}k and E\_{}k]{Specification of $\boldsymbol{\mu_k}$ and $\boldsymbol{E_k}$}\label{tauketc}

Since $M$ is non-compact with compact boundary, there exists a smooth, proper map $\rho\colon M\rightarrow (-\infty,0]$ (see \cite[Lemma 2.1]{Zhu2020}) such that
\[
	\rho^{-1}(0)=\partial M,\qquad |\ed\rho|_g<1.
\]
Fix a smooth function $\eta\in C^\infty((-\infty,0])$ satisfying
\[
	\eta(t)=0 \quad\mbox{for any } t\le -1 \qquad \mbox{and}\qquad \eta (0)=2;
\]
define $\tau_k$ by
$ 3\coth(2\tau_k) = 3+ k^{-1}$,
and then define $\hat\mu_{k}\colon (-\tau_k,0]\rightarrow \mathbb{R}$ by
\[
	\hat\mu_{k}(t)=3\coth\bigl(2(t+\tau_k)\bigr)-k^{-1}\eta(t).
\]
Thus, $\{\tau_k\}_{k = 1}^\infty$ is increasing and tends to infinity, and
\[
	\hat\mu_{k}(-\tau_k)=+\infty\qquad\mbox{and}\qquad \hat\mu_{k}(0)=3-k^{-1}.
\]
Now, choose $a_k$, regular values of
$\rho$, such that $\tau_k\le a_k < \min\{\tau_{k+1}, \tau_k + 1\}$, and then define $E_{k}:=\rho^{-1}\bigl([-a_{k},0]\bigr)\subset M$.
Denote $\partial_{k}^- :=\rho^{-1}(-a_{k})$,
which are smooth hypersurfaces of $M$. This makes $(E_k, g|_{E_k}; \partial_k^-,\partial M)$
a Riemannian band.
Finally, let $\rho_{k}:=(\tau_k/a_k)\rho$, and define
 \[
 	\mu_k := \hat\mu_k\circ\rho_k.
\] By this arrangement, $\mu_k|_{\partial_k^-} = \infty$.

\subsection[mu\_{}k-bubbles in E\_{}k]{$\boldsymbol{\mu_k}$-bubbles in $\boldsymbol{E_k}$}
\label{muk_in_Ek}
For each fixed $k$, consider $(E_k, g|_{E_k}; \partial_k^-, \partial M)$.
Note that $H_{\partial M}\ge 3$; by construction, $\mu_k$ satisfies
the barrier condition (see Definition~\ref{barrierdef}).
By Fact~\ref{mubbExist}, a smooth $\mu_k$-bubble $\Omega_k$ exists.
Define $\Sigma_k := \partial\Omega_k \setminus \partial_k^-$, which is smooth, closed, and separates $\partial_k^-$ from $\partial M$.

The following lemma shows that all $\Sigma_k$ must meet a
fixed compact subset of $M$.

\begin{Lemma}\label{intersectionLemma}
	Let $\Kcal:= \{x\in M\colon \dist_g(x,\partial M)\le 10\}$.
	Then $\Sigma_k\cap \Kcal\ne \varnothing$.
\end{Lemma}
\begin{proof}
Suppose on the contrary that $\Sigma_{k}\cap \Kcal = \varnothing$. This implies that
$\eta\circ\rho_k = 0$ on $\Sigma_k$.
Moreover, by assumption, $R_g\ge -12$, and by construction, $|\ed\rho_k|_g < 1$. Thus, we have (see \eqref{RplusDef})
\[
R_{+}^{\mu_{k}}> -12+\frac{4}{3}[3\coth(2(\rho_{k}+\tau_k))]^{2}-12[\sinh(2(\rho_{k}+\tau_k)]^{-2}=0\qquad  \text{on }  \Sigma_k.
\]
This, along with Fact~\ref{confPSC}, implies that $\Sigma_k$ admits a PSC metric;
since $\Sigma_k$ is separating, we get a contradiction, by Lemma~\ref{4dNSepLemma}.
\end{proof}

\subsection[Convergence of Sigma\_{}k]{Convergence of $\boldsymbol{\Sigma_k}$}
\label{Sigmak_limits}
By using \cite[Theorem 3.6]{ZZ20}, one can show that the
 second fundamental form $\II_{\Sigma_k}$
is uniformly bounded within any compact subset of $M$.
Thus, by Lemma~\ref{intersectionLemma}, $\Sigma_k$ subconverges to a smooth hypersurface $\Sigma$ in $M$ (for convenience, denote the subsequence by the same symbol $\Sigma_k$). Within
compact subsets of $M$, the convergence is uniform and has multiplicity one;
moreover, $\Sigma$ bounds a `minimizing $3$-bubble' for which minimality is interpreted
with respect to compactly supported
perturbations {(cf. \cite[Lemma 4.10]{HLS22})}. Depending on whether $\Sigma$ is compact, we consider the
two cases below.

\emph{Case $1$: $\Sigma$ is compact}. By minimality, we have (see Fact~\ref{mubbVarProp})
\begin{equation}\label{stability_cptCase}
	H_{\Sigma}=3\qquad
	\mbox{and}\qquad L_{\Sigma}=-\Delta_{\Sigma}+\frac{1}{2}\big(R_{\Sigma}-R_{+}^{3}\big)\geq 0.
\end{equation}
Since $R_{+}^{3}=R_{g}+12\geq 0$, \eqref{stability_cptCase} implies that
$-\Delta_{\Sigma}+\frac{1}{2}R_{\Sigma}\ge 0$; thus, there exists
a smooth function $u>0$ defined on $\Sigma$ and a constant $\lambda\ge0$
such that
\begin{equation}\label{eigenfcn}
	\left(-\Delta_{\Sigma}+\frac{1}{2}R_{\Sigma}\right)u=\lambda u.
\end{equation}	
Define $\tilde{g}_{\Sigma}=ug_{\Sigma}$ where
$g_\Sigma$ is the metric on $\Sigma$ induced by $g$. We have
\begin{equation}\label{tildegScal}
R_{\tilde{g}_{\Sigma}}=u^{-1}\left(R_{{\Sigma}}+\frac{3}{2}\left|\frac{\nabla u}{u}\right|^{2}-2\frac{\Delta u}{u}\right)=u^{-1}\left(2\lambda+\frac{3}{2}\left|\frac{\nabla u}{u}\right|^{2}\right)\geq 0.
\end{equation}

Since each $\Sigma_k$ is separating, so is $\Sigma$. By Lemma~\ref{4dNSepLemma},
 $\Sigma$ admits no PSC metric; then by~\eqref{tildegScal} and
 the trichotomy theorem of Kazdan and Warner,
 $R_{\tilde g_\Sigma}= 0$. Thus,  $\lambda$ must vanish, and~$u$~must be a constant; \eqref{eigenfcn} in turn implies that $R_{\Sigma}=0$. %and  \eqref{stability_cptCase} enforces that $R_g = -12$ on $\Sigma$.
Then by Bourguignon's theorem (see \cite[Lemma 5.2]{KW75}),
$g_\Sigma$ is Ricci-flat, which must be flat since $\dim \Sigma=3$.

Now we prove that a neighborhood of $\Sigma$ splits.
When $\Sigma\cap\partial M = \varnothing$, since $\Sigma$ is the boundary of minimizing $3$-bubble,
\cite[Theorem 2.3]{ACG08} implies that  there exists an open neighborhood of $\Sigma$ that is isometric to a warped product $\bigl((-\epsilon,\epsilon)\times\Sigma, \ed t^{2}+{\rm e}^{2t}g_{\Sigma}\bigr)$, where $t$ is the coordinate on $(-\epsilon,\epsilon)$ and $\Sigma$ corresponds to $t = 0$.  When $\Sigma\cap\partial M\ne \varnothing$, we must have $\Sigma = \partial M$, by the maximum principle.
In this case, the proof of \cite[Theorem 2.3]{ACG08} still applies and gives an open neighborhood of $\Sigma$ that is isometric to a warped product $\bigl((-\epsilon,0]\times\Sigma, \ed t^{2}+{\rm e}^{2t}g_{\Sigma}\bigr)$.

Thus, a neighborhood of $\Sigma$ is foliated by the $t$-level sets. Note that moving along the foliation leaves the energy functional invariant; thus, each $t$-slice also bounds a minimizing $3$-bubble, to which the same analysis above applies.

This implies that a maximal neighborhood $\Ucal$ of $\Sigma$ on which the metric splits as
\[
	\bigl(I\times\Sigma, \ed t^{2}+{\rm e}^{2t}g_{\Sigma}\bigr)
\]  must be both
open and closed in $M$. By connectedness, $\Ucal = M$, and $I$ must be of the form $(-\infty, c]$.
This achieves the desired splitting.

\emph{Case $2$: $\Sigma$ is non-compact}. By finding a contradiction, we prove that this case does not occur. The argument largely follows the proof of \cite[Theorem 1.1]{Zhu2020}, so we only sketch the steps.

Let
\[
	\left(M_k, g_k\right)=\bigl(\Sigma_{k}\times \S^{1}, g_{\Sigma_{k}}+u_{k}^{2}dt^{2}\bigr),
\]	
where $u_{k}$ is the first eigenfunction of $L_{\Sigma_{k}}$; that is, $L_{\Sigma_{k}}u_{k}=\lambda_{k}u_k$ with $\lambda_{k}\geq 0$.
Since $\dim \Sigma_k = 3$,  \cite[Corollary 1.10]{CRZ2022} implies that $M_k$ admits no PSC metric.

Now
\begin{equation}\label{Rgk}
R_{g_{k}}=R_{g_{\Sigma_{k}}}-2\frac{\Delta_{g_{\Sigma_{k}}}u_{k}}{u_{k}}=R_{+}^{\mu_{k}}+2\lambda_{k}.
\end{equation}
By construction, $R_{+}^{\mu_{k}} \ge 0$ outside $\Kcal$, and
there exist $\delta_k>0$, satisfying $\lim\delta_k= 0$, such that $R_{+}^{\mu_{k}}\geq -\delta_{k}$ on $M$.
Since $R_{g_k}$ cannot be positive and $\lambda_k\ge 0$, by \eqref{Rgk}, we must have $\lim \lambda_{k}= 0$.

Next, choose $q_{k}\in \Sigma_{k}\cap \Kcal$ so that $\lim q_{k} = q\in\Sigma$, and let $p_{k}=(q_{k},t_{0})\in \Sigma_k\times \S^1$ and $p=(q,t_{0})\in \widetilde M=\Sigma\times \S^{1}$.
Normalize $u_{k}$ such that $u_{k}(q_{k})=1$. By the Harnack
inequality, $u_k$ converges smoothly to a positive function $u$ on $\Sigma$ with $u(q)= 1$.
Thus, $(M_k, g_k)$ converges in the pointed smooth topology to $\bigl(\widetilde M, \widetilde g\bigr)$, where $\widetilde g = g_{\Sigma}+ u^2\ed t^2$.

Now one can follow the proof\,%
\footnote{The proof of \cite[Proposition 3.2]{Zhu2020} only relies on $\widetilde M$ admitting no PSC metric and the properties of $R^{\mu_k}_+$ mentioned above.}
of \cite[Proposition 3.2]{Zhu2020} to show
that $\Ric_{\widetilde g} = 0$,
and then follow the proof\,%
\footnote{In particular, the boundedness of  $\area(\Sigma)$ follows from
$
	\Acal^{\mu_k}_{\refOmega}(\Omega_{k})\leq \Acal^{\mu_k}_{\refOmega}(E_{k})
$ and $\mu_{k}>0$.}
of \cite[Theorem 1.1]{Zhu2020} to show that $u$ is constant, which
implies $\Ric_{g_{\Sigma}}= 0$.

In summary, $(\Sigma, g_\Sigma)$ is complete, non-compact, Ricci-flat, and with finite area; this contradicts \cite[p.~25, Theorem 4.1]{Yau94}.
\qed

\begin{Remark}\label{warpProdRmk}
The PSC obstruction, provided by \cite[Corollary 1.10]{CRZ2022}, for manifolds of the form $\Sigma\times \S^1$ only works when $\dim\Sigma \ne 4$. On the other hand, if $\Sigma$
$(2\le\dim\Sigma\le 6)$ is closed, orientable, and if it admits a map of nonzero degree to some
$\Sigma'\in \Ccal_{\deg}$, then by a similar argument as \cite[Theorem 1.1]{CLSZ2021}, one can show that $\Sigma\times \S^{1}$ admits no PSC metric.
 \end{Remark}

\begin{proof}[Proof of Theorem~\ref{alhpmt1}]
The inequality
$
	\inf_{\partial M}H \le n-1$
follows directly from Proposition~\ref{incprssArg}.

To prove the second part of the theorem, first obtain a covering $\bigl(\hat M, \hat g\bigr)$ of
$(M,g)$ as in the proof of Proposition~\ref{incprssArg}, and then apply (essentially)
the same proof of Theorem~\ref{alhpmt2} to $\big(\hat M,\hat g\big)$; to assist the reader,
we list a few points that may need attention.
\begin{itemize}\itemsep=0pt
	\item{$\partial \hat M$ may not be connected, but Riemannian bands can still be constructed in a similar manner as in
			the proof of Proposition~\ref{incprssArg}. To avoid clash of symbols, denote $S:=\partial M$ and let $\hat S$ be
			a fixed lifting of $S$ in $\hat M$. Thus $\partial_+ = \hat S$ and $\partial_\star\subset\partial_-$;
			$\mu_k >0$ can be defined such that $\mu_k|_{\partial_\star} = \infty$ and $\mu_k|_{\hat S} = (n-1) - 1/k$;
			on $\partial_-\setminus\partial_\star$ (if nonempty) we have $H \ge n-1$; one can check that
			the barrier condition is satisfied, and the $\Sigma_k$s exist; restricting $j\colon \hat M\rightarrow S$ to $\Sigma_k$
			yields a map $\Sigma_k\rightarrow\hat S$ of nonzero degree.}
	\item{An adapted version of Lemma~\ref{intersectionLemma} holds; in the proof, invoke Lemma~\ref{C_degNSepLemma}
	instead of Lemma~\ref{4dNSepLemma}. It follows that $\Sigma_k$ converges to some $\Sigma$.}		
	\item{When $\Sigma$ is compact, the corresponding part in Section~\ref{Sigmak_limits} applies, apart from dimensional adjustments
		and the fact that Ricci-flatness may no longer imply flatness.}
	\item{When $\Sigma$ is non-compact, we need to argue, without relying
	 on \cite[Corollary 1.10]{CRZ2022},  that $M_k = \Sigma_k\times\S^1$ admits no PSC metric, and this is already addressed by Remark~\ref{warpProdRmk}.}
\end{itemize}
The consequence is that
 $\big(\hat M, \hat g\big)$ is of the form
 \[
 	\bigl((-\infty, 0]\times \Sigma, \ed t^2+{\rm e}^{2t}g_\Sigma\bigr),
\]
where
  $g_\Sigma$ is Ricci-flat. In particular, the covering $\hat M\rightarrow M$ is $1$-fold and hence an isometry.
  Since $\Sigma = \partial M$ is assumed to be aspherical,
  $g_\Sigma$ must be flat, which can be seen by applying the Cheeger--Gromoll splitting theorem
  to the universal cover; for details, see the beginning paragraph of \cite[Section~6]{CL20}.
\end{proof}

\appendix

\section[mu-bubbles]{$\boldsymbol{\mu}$-bubbles}\label{mubbSec}

This section collects some `definitions' and `facts' concerning
the $\mu$-bubble technique, about which we make no claim to originality.
For detailed expositions and proofs,
the reader may consult \cite{CRZ2022,CL20,ZZ20, Zhu21warp} and \cite[Section~5]{Gr2021}.
This section also includes three supplementary `lemmas'.\looseness=-1

A common setting for $\mu$-bubbles is a \emph{Riemannian band}, namely
a compact, connected Riemannian manifold $\big(M^n,g\big)$ whose (nonempty) boundary
is expressed as a disjoint union $\partial M = \partial_-\sqcup \partial_+$,
where each of $\partial_\pm$ is a smooth, closed and possibly disconnected
$(n-1)$-manifold.

Given a Riemannian band $\big(M^n,g;\partial_-,\partial_+\big)$
and a function $\mu\in C^\infty\big(\mathring{M}\big)$,
consider the following variational problem: Let $\Omega_0$ be a smooth open neighborhood of $\partial_-$; among all Caccioppoli sets
$\Omega\subset M$ that satisfy $\partial_-\subset \Omega$ and $\Omega\Delta\Omega_0 \Subset \mathring{M}$, seek a \emph{minimizer} of the functional
\begin{equation*}%\label{braneAct}
	\Acal_{\refOmega}^{\mu}(\Omega)=\Hcal^{n-1}(\partial \Omega)-\Hcal^{n-1}(\partial\refOmega)-\int_{M}(\chi_{\Omega} - \chi_{\refOmega})\mu\, {\rm d}\Hcal^n,
\end{equation*}
where $\Hcal^k$ is the induced $k$-dimensional Hausdorff measure, and
$\chi_\Omega, \chi_{\refOmega}$ are characteristic functions.
Such a minimizer is called a \emph{$\mu$-bubble}.

Existence and regularity of $\mu$-bubbles are well-established when $\mu$
satisfies the following `barrier condition'.
\begin{Definition}\label{barrierdef}
	Let $\big(M^n,g;\partial_-,\partial_+\big)$ be a Riemannian band.
	A function $\mu\in C^\infty\big(\mathring{M}\big)$ is said to satisfy the
	\emph{barrier condition} if, for each connected component $S\subset\partial_+$ (resp., $S\subset \partial_-$),
	\begin{itemize}\itemsep=0pt
	\item either $\mu$ smoothly extends to $S$ and satisfies
					$H_S > \mu|_S$	(resp., $H_S > -\mu|_S$),
	where $H_S$ is
	the mean curvature of $S$ with respect to the outward normal;
	\item or $\mu\to -\infty$  (resp., $\mu\to +\infty$) towards $S$.
	\end{itemize}
\end{Definition}

\begin{Fact}\label{mubbExist}
	For $3\le n\le 7$, if $\mu\in C^\infty\big(\mathring{M}\big)$ satisfies the barrier condition, then
	there exists a~smooth $\mu$-bubble $\Omega$.
	In particular, $\partial\Omega \setminus \partial_-$ is homologous to
	$\partial_+$ and is separating $($see Definition~{\rm \ref{noSepDef}} below$)$.
\end{Fact}

Also well-known are the following variational properties. To fix notation,
let $\Sigma$ denote the hypersurface $\partial\Omega\setminus \partial_-$
with outward unit normal $\nu$; let $R_\Sigma$ and $\Delta_\Sigma$ be, respectively, the scalar
curvature and the Laplacian along $\Sigma$ (with the induced metric);
 let $H_\Sigma$ and $\II$ be, respectively, the mean curvature and
the second fundamental form of $\Sigma$, computed with respect to $\nu$;
define the operators
\begin{gather*}%\label{JDef}
	J_\Sigma = -\Delta_\Sigma + \frac{1}{2}\bigl(R_\Sigma - R_g - \mu^2 - |\II|^2\bigr) - \nu(\mu)
\\ %\label{LDef}
	L_\Sigma =  - \Delta_\Sigma + \frac{1}{2}\bigl(R_\Sigma - R^\mu_+\bigr),
\end{gather*}
where
\begin{equation}\label{RplusDef}
	R^\mu_+ = R_g + \frac{n}{n-1}\mu^2 - 2|\ed\mu|_g.
\end{equation}
\begin{Fact}\label{mubbVarProp}
	Suppose that $\Omega$ is a smooth $\mu$-bubble.
	We have
	\begin{enumerate}\itemsep=0pt
		\item[$(a)$] $H_\Sigma = \mu|_\Sigma$;
		\item[$(b)$] $L_\Sigma \ge J_\Sigma\ge 0$.
	\end{enumerate}
\end{Fact}

The semi-positivity of $L_\Sigma$ has several applications, and we shall list a few.
To start with, let $u > 0$ be an eigenfunction associated to the first eigenvalue
$\lambda\ge 0$ of $L_\Sigma$. Consider the warped-product metric
$\hat h :=g_\Sigma + u^2\ed \theta^2$ defined on $\hat \Sigma := \Sigma\times \S^1$,
where $\theta\in \S^1$.
\begin{Fact}\label{warpPSC}
	Suppose that $\Omega$ is a smooth $\mu$-bubble.
	The scalar curvature of $\big(\hat\Sigma,\hat h\big)$
	is
	\begin{equation*}
		R_{\hat h} = R_\Sigma - 2u^{-1}\Delta_\Sigma u = R^\mu_+ +  2\lambda.
	\end{equation*}
	In particular, if $R^\mu_+ >0$ on $\Sigma$, then $\Sigma\times \S^1$
	admits a~PSC metric.
\end{Fact}

Alternatively, one can compare $L_\Sigma$ with the conformal Laplacian on $\Sigma$ and obtain the following.
\begin{Fact}\label{confPSC}
	For $n\ge 3$, suppose that $\Omega$ is a smooth $\mu$-bubble on which
	$R^\mu_+ >0$. Then $\Sigma$ admits a~PSC metric.
\end{Fact}

With additional topological assumptions on $M$, Fact~\ref{confPSC} can be used
to prove width estimates for $(M,g)$. To be precise, we start by recalling the following notion (cf. \cite[Property~A]{CRZ2022}).
\begin{Definition}\label{noSepDef}
	Given a (topological) band $\big(M^n;\partial_-, \partial_+\big)$, we say that
	a~(closed) hypersurface~$\Scal$ in~$M$ is
	\emph{separating}, if all paths connecting
	$\partial_-$ and $\partial_+$ must intersect $\Scal$.
	A band is said to satisfy
	the \emph{\Nsep property} if no separating hypersurface admits a PSC metric.
\end{Definition}

\begin{Remark}\label{sepRmk}
	If $\Scal\subset M^n$ is a separating hypersurface, then there exists a
	minimal list of connected components $S_i$ $(i = 1,\ldots, k)$ of $\Scal$
	such that their union $\Scal'$ remains separating.
	For details, see \cite[Lemma 2.2]{CRZ2022}.
	Using intersection theory,
	one can show that
	$[\Scal']\ne 0\in H_{n-1}(M;\Z)$. Moreover, with suitable orientation, $\Scal'$ is homologous to $\partial_{+}$ in $M$.
\end{Remark}

\begin{Lemma}\label{C_degNSepLemma}
	Let $\big(M^n, g; \partial_-, \partial_+\big)$ be a
	Riemannian band, and let $\iota\colon \partial_+\hookrightarrow M$ be the inclusion map.
Suppose that $\partial_+ \in \mathcal{C}_{\deg}$ $($see Definition~{\rm \ref{cdegdef})} and that there exists a continuous map ${f\colon M\rightarrow \partial_+}$ such that $f\circ\iota$ is homotopic to  $\id_{\partial_+}$. Then
$(M,g)$ satisfies the \Nsep property.
\end{Lemma}
\begin{proof}Suppose that $\Scal$ is a separating hypersurface in $M$,  and  let $\Scal'$ be as in Remark \ref{sepRmk}; in particular,
 $\Scal'$  is homologous to $\partial_+$ in $M$. Now since $f\circ\iota$ is homotopic to $\id_{\partial_+}$,  it is easy to see that the restriction $f|_{\Scal'}\colon \Scal'\rightarrow \partial M$ has degree $1$. Since $\partial_+\in \mathcal{C}_{\deg}$, $\Scal'$ admits no PSC metric.
\end{proof}

The next lemma is a variant of Gromov's band-width estimate \cite[Section 5.3]{Gr2021}.

\begin{Lemma}\label{bandWidth}
	For $3\le n\le 7$, let $\big(M^n,g;\partial_-,\partial_+\big)$ be a Riemannian band that satisfies
	the \Nsep property, and let
	$\partial_\star\subset \partial_-$ be a compact subset without boundary.
	Suppose that
	\begin{enumerate}\itemsep=0pt
		\item[$(a)$] $R_g \ge -n(n-1)$;
		\item[$(b)$] $H_{\partial_-\setminus\partial_\star} \ge -(n-1)$;
		\item[$(c)$] $H_{\partial_+} \ge (n-1)( 1+ \delta)$ for some constant $\delta>0$.
	\end{enumerate}
	Then there exists a constant $T_\delta>0$, depending only on $\delta$, such that
	\begin{equation*}%\label{mixedWthEst}
		\dist_g(\partial_\star,\partial_+)
		\le T_\delta.
	\end{equation*}	
\end{Lemma}
\begin{proof}
	Set  $\epsilon = \delta/3$,
	and define
	$C_\delta, T_\delta>0$ by
	\[
		  \coth(C_\delta/2) = \frac{1+\delta/2}{1+\epsilon}
		  \qquad\mbox{and}
		  \qquad T_\delta = \frac{C_\delta}{n(1+\epsilon)}.
	\]
	For the sake of deriving a contradiction, suppose that
	$\dist_g(\partial_\star, \partial_+) > T_\delta$.
	By the proof of \cite[Lemma~4.1]{Zhu21warp}, there exists a smooth, proper function
	$\rho\colon M\rightarrow [-T_\delta, 0]$
	such that
	\begin{equation}\label{rhofcnprop}
		\rho^{-1}(-T_\delta) = \partial_\star, \qquad \rho^{-1}(0) = \partial_+,
		\qquad\mbox{and}\qquad |\ed\rho|_g < 1.
	\end{equation}
	Now consider the function
	\[
		h(t)= (n-1)(1+\epsilon)\coth\left(\frac{n(1+\epsilon)t + C_\delta}{2}\right),
		\qquad t\in (-T_\delta, 0].
	\]
	By construction, $h$ is decreasing, strictly greater than $n-1$, and satisfies
	\begin{equation}\label{hfcnprop}
		h(0) < H_{\partial_+}, \qquad \lim_{t\rightarrow -T_\delta} h(t) = \infty,
		\qquad \frac{n}{n-1}h(t)^2 + 2h'(t) \equiv n(n-1)(1+\epsilon)^2.
	\end{equation}
	Combining \eqref{rhofcnprop}, \eqref{hfcnprop}, and the assumptions $(a)$, $(b)$, $(c)$,
	one can easily check that the
	function $\mu:=h\circ\rho$, defined on $M\setminus \partial_\star$,
	satisfies both
	the barrier condition and the inequality $R^\mu_+ >0$.
	By Facts~\ref{mubbExist} and \ref{confPSC}, there exists a separating hypersurface $\Sigma$ in $(M;\partial_-,\partial_+)$ that admits a PSC metric. This contradicts the \Nsep hypothesis.
\end{proof}

\begin{Remark}\label{bandNonexistence}
	We mention two variants of Lemma~\ref{bandWidth}, both of which can be obtained
	by slightly modifying the proof above.
	$(A)$ For $3\le n\le 7$,
	no Riemannian band can simultaneously satisfy the \Nsep property and
	the conditions
	$R_g \ge -n(n-1)$, $H_{\partial_-} \ge -(n-1)$ and $H_{\partial_+}> n-1$.
	$(B)$ For $3\le n\le 7$, let $(M^n,g)$ be a complete, non-compact Riemannian
	manifold with compact boundary $\partial M$. Suppose that
	$M$ satisfies the \Nsep property (see below); then $(M,g)$ cannot satisfy
	the conditions $R_g \ge -n(n-1)$ and $H_{\partial M} > n-1$ simultaneously.
\end{Remark}

The concept of separating hypersurfaces can also be defined for complete, non-compact Riemannian manifolds $(M,g)$ with compact boundary---just require that $\Scal$ intersects with all paths connecting $\partial M$ and infinity. The
\Nsep property can be extended to such manifolds.

\begin{Lemma}\label{4dNSepLemma}
	Let $\big(M^4, g\big)$ be a complete, non-compact Riemannian $4$-manifold with compact $($nonempty$)$ boundary $\partial M$.
Suppose that the homotopy groups $\pi_2(M) = \pi_3(M) = 0$. Then~$(M,g)$ satisfies the \Nsep property.
\end{Lemma}
\begin{proof}
	Suppose that $\Scal\subset M$ is a (closed) separating hypersurface that
	admits a PSC metric, and let $\Scal'\subset \Scal$ be as indicated
	in Remark~\ref{sepRmk}. In particular, $\Scal'$ admits a PSC metric, and
	 $[\Scal']\ne 0\in H_3(M,\Z)$.
	Since $\pi_2(M)$ is trivial, the topological classification of
	closed 3-manifolds admitting
	a PSC metric implies that
	$\Scal'$ is homologous to a spherical class in $H_{3}(M,\Z)$
	(see \cite[p. 112]{WangJthesis}). Since $\pi_3(M)$ is also trivial, this violates
	Lemma~\ref{topLemma} below.
\end{proof}

\section{Topological lemmas}\label{TopologicalLemmas}

\begin{Lemma}\label{topLemma}
	Let $M$ be a non-compact $4$-manifold satisfying $\pi_3(M) = 0$. Then
	$H_3(M,\Z)$ contains no nontrivial spherical class $($i.e., classes of the form
	$\bigl[\S^3/\Gamma\bigr])$.
\end{Lemma}
\begin{proof}
	Let $[\beta]$ denote the fundamental class of $\S^3/\Gamma$
	where $\Gamma$ is a discrete subgroup of $O(4)$. Let $i\colon \S^3/\Gamma\rightarrow M$ be a continuous map. The goal is
	to prove that $i_*[\beta] = 0\in H_3(M,\Z)$. Now let $[\alpha]$
	be the fundamental class of $\S^3$. The composition
	$\S^3\xrightarrow{\pi}\S^3/\Gamma\xrightarrow{i} M$ induces a
	map at the level of $H_3(\cdot, \Z)$, such that
	$[\alpha]\xrightarrow{\pi_*} d[\beta] \xrightarrow{i_*}di_*[\beta]$
	where $d$ is the degree of $\pi$. Since $\pi_3(M) = 0$, Hurewicz homomorphism
	implies that \[
		di_*[\beta] = (i\circ\pi)_*[\alpha] = h([i\circ\pi])= 0\in H_3(M,\Z),
	\]
	where $h\colon \pi_3(M)\rightarrow H_3(M,\Z)$ is the Hurewicz map. Thus,
	in order to show that $i_*[\beta] = 0$, it suffices to show that
	$H_3(M,\Z)$ is torsion free, and this follows from $M$ being non-compact
	(see \cite[Corollary~7.12]{bredon}).
\end{proof}

\begin{Lemma}\label{incprLemma}
	For $1\le k\le n-2$, let $M_1 = \big(\R\times \T^{n-k-1} - \Bb\big)\times \T^k$, where $\Bb$ is an embedded $(n-k)$-ball in $\R\times \T^{n-k-1}$.
	Let $M_2$ be a smooth, possibly non-compact, manifold with boundary~$\partial M_2$.
	Suppose that  $\Phi\colon  \partial M_1 \rightarrow \partial M_2$ is a diffeomorphism, and let
	$M:= M_1\sqcup_\Phi M_2$ be the manifold obtained by identifying $\partial M_1, \partial M_2$ via $\Phi$.
	Let $t\in \R$ be such that $\{t\}\times \T^{n-k-1}$ is disjoint from~$\Bb$.
	Then the hypersurface $\Sigma = \{t\}\times \T^{n-1}$ is incompressible in $M$
	if and only if the $\T^k$-factor\,\footnote{Note that $\partial M_2$ has the product structure $\S^{n-k-1}\times\T^k$ induced by $\Phi$.}
	of~$\partial M_2$ is incompressible in $M$.
\end{Lemma}
\begin{proof}
	In $M$, the $\T^k$-factor of $\Sigma$ is homotopic to that of $\partial M_1$ and hence to that of $\partial M_2$. Thus, $(\Rightarrow)$
	is clear.
	
	For $(\Leftarrow)$, we prove its contrapositive. Suppose that
	 $L\subset \Sigma$ is a non-contractible loop that is contractible in $M$. Write
	 \[
	 	[L] = (m_i\alpha_i, n_j\beta_j) \in\pi_1\big(\T^{n-k-1}\big)\times \pi_1\big(\T^{k}\big)\cong \pi_1(\Sigma),
	 \]
	 where $\alpha_i$ generates the fundamental group of the $i$-th $\S^1$-factor in $\T^{n-k-1}$ and $m_i\in \Z$,
	 similarly for $\beta_j$ and $n_j$.
	 Let us write $\hat\alpha_i$, $\hat \beta_j$ for the corresponding elements in the homology class $H_1(\Sigma;\Z)$.
	
	 It will be convenient to view the $\R$-factor in $\R\times \T^{n-k-1}$ as $\S^1$ minus a point, and to view
	 $M$ as a subset of $\hat M:= \hat M_1\sqcup_\Phi M_2$, where  $\hat M_1 := \bigl(\S^1\times \T^{n-k-1} - \Bb\bigr)\times \T^k$.
		
	 Let $\iota\colon \Sigma\hookrightarrow \hat M$ be the inclusion map.
	 For $1\le i\le n-k-1$, let $\theta_i$ be the coordinate on the $i$-th $\S^1$-factor of $\T^{n-k-1}$. By construction, there exists
	 $t_i\in \S^1$ such that $\theta_i = t_i$ defines a~hypersurface $S_i$ in $\hat M$ that is `dual' to $\iota_*\hat\alpha_i$, in the sense that
	 the intersection products
	 \[
	 	[S_i]\cdot \iota_*\hat \alpha_i = 1\qquad \mbox{and}\qquad  [S_i]\cdot \iota_*\hat\alpha_{i'} =[S_i]\cdot \iota_*\hat\beta_j = 0,
		\qquad i'\ne i.
	\]
	Since $L$ is contractible in $M\subset \hat M$, we have
	\[
		\sum_i m_i \iota_*\hat \alpha_i + \sum_j n_j\iota_*\hat \beta_j = 0 \in H_1(\hat M;\Z);
	\]
	by taking intersection products with $[S_i]$, we see that $m_i = 0$ for all $i = 1,\ldots, n-k-1$,
	so $L$ is homotopic to a loop in the $\T^{k}$-factor of $\Sigma$. Thus, the $\T^k$-factor of $\Sigma$
	is not incompressible in~$M$. By homotopy, the same is true for the $\T^k$-factor of $\partial M_2$.
	This completes the proof.
\end{proof}

\subsection*{Acknowledgements}

We thank Shihang He for kindly sharing his proof of Lemma~\ref{topLemma}. We also thank the anonymous referees for carefully reading the manuscript and offering suggestions, which has led to improved exposition and a more direct argument now included in Remark~\ref{vanKampenRmk}.
Research leading to this work was supported by the National Key R\&D Program of China Grant 2020YFA0712800 (T.~Hao, P.~Liu and Y.~Shi) and the China Postdoctoral Science Foundation Grant 2021TQ0014 (Y.~Hu).

\pdfbookmark[1]{References}{ref}
\LastPageEnding

\end{document}